\documentclass[11pt,twoside]{article}
\usepackage{amsmath, amsthm, amscd, amsfonts, amssymb, graphicx, color}

\date{}

\begin{document}

\centerline{}

\centerline {\Large{\bf Controlled frames in \,$n$-Hilbert spaces and their tensor products }}

\newcommand{\mvec}[1]{\mbox{\bfseries\itshape #1}}
\centerline{}
\centerline{\textbf{Prasenjit Ghosh}}
\centerline{Department of Pure Mathematics, University of Calcutta,}
\centerline{35, Ballygunge Circular Road, Kolkata, 700019, West Bengal, India}
\centerline{e-mail: prasenjitpuremath@gmail.com}
\centerline{}
\centerline{\textbf{T. K. Samanta}}
\centerline{Department of Mathematics, Uluberia College,}
\centerline{Uluberia, Howrah, 711315,  West Bengal, India}
\centerline{e-mail: mumpu$_{-}$tapas5@yahoo.co.in}

\newtheorem{Theorem}{\quad Theorem}[section]

\newtheorem{definition}[Theorem]{\quad Definition}

\newtheorem{theorem}[Theorem]{\quad Theorem}

\newtheorem{remark}[Theorem]{\quad Remark}

\newtheorem{corollary}[Theorem]{\quad Corollary}

\newtheorem{note}[Theorem]{\quad Note}

\newtheorem{lemma}[Theorem]{\quad Lemma}

\newtheorem{example}[Theorem]{\quad Example}

\newtheorem{result}[Theorem]{\quad Result}
\newtheorem{conclusion}[Theorem]{\quad Conclusion}

\newtheorem{proposition}[Theorem]{\quad Proposition}

\begin{abstract}
\textbf{\emph{The concepts of controlled frames and it's dual in n-Hilbert spaces and their tensor products have been introduced and then some of their characterizations are given.\,We further study the relationship between controlled frame and bounded linear operator in tensor product of n-Hilbert spaces.\,At the end, the direct sum of controlled frames in \,$n$-Hilbert space is being considered.}}
\end{abstract}
{\bf Keywords:} \emph{n-normed space, n-inner product space, tensor product of Hilbert spaces,\\ \smallskip\hspace{2.5cm} frame, dual frame, controlled frame.}


\section{Introduction}

\smallskip\hspace{.6 cm} In the study of vector spaces, the main characterestic of a basis is that every elements can be represented as a superpostion of the elements in the basis.\;A frame is also sequence of elements in a Hilbert space, which allows every elements to be written as a linear combinations of the elements in the frame.\,However, the corresponding coefficients are not necessarily unique.\;So, frame can be considered as a generalizations of a basis.\;In fact, frames play important role in theoretical research of wavelet analysis, signal denoising, feature extraction, robust signal processing etc.\;In 1946, D.\,Gabor \cite{Gabor} first initiated a technique for rebuilding signals using a family of elementary signals.\;In 1952, Duffin and Schaeffer abstracted Gabor's method to define frame for Hilbert space in their fundamental paper \cite{Duffin}.\,Later on, in 1986, as the wavelet era began, Daubechies et. al \cite{Daubechies} observed that frames can be used to find series expansions of functions in \,$L^{\,2}\,(\,\mathbb{R}\,)$\, which are very similar to the expansions using orthonormal bases.

\,Controlled frame is one of the newest generalization of frame.\,I. Bogdanova et al.\,\cite{I} introduced controlled frame for spherical wavelets to get numerically more efficient approximation algorithm.\,Thereafter, P. Balaz \cite{B} developed weighted and controlled frame in Hilbert space.\,S.\,Rabinson \cite{S} presented the basic concepts of tensor product of Hilbert spaces.\,The tensor product of Hilbert spaces \,$X$\, and \,$Y$\, is a certain linear space of operators which was represented by Folland in \cite{Folland}, Kadison and Ringrose in \cite{Kadison}.

\;The concept of \,$2$-inner product space was first introduced by Diminnie et. al \cite{Diminnie} in 1970's.\;In 1989, A.\,Misiak \cite{Misiak} developed the generalization of a \,$2$-inner product space for \,$n \,\geq\, 2$.

In this paper, our focus is to study controlled frames associated to \,$\left(\,a_{\,2},\, \cdots,\, a_{\,n}\,\right)$\, in \,$n$-Hilbert spaces and their tensor products.\,We will see that any controlled frame associated to \,$\left(\,a_{\,2},\, \cdots,\, a_{\,n}\,\right)$\, is a frame associated to \,$\left(\,a_{\,2},\, \cdots,\, a_{\,n}\,\right)$\, in \,$n$-Hilbert space and the converse part is also true under some conditions.\,In tensor product of \,$n$-Hilbert spaces, we shall established that an image of a controlled frame under a bounded linear operator will be a controlled frame if and only if the operator is invertible.\,Dual controlled frame in tensor product of \,$n$-Hilbert spaces is discussed and finally, we shall established that the direct sum of controlled frames associated to \,$\left(\,a_{\,2},\, \cdots,\, a_{\,n}\,\right)$\, is again a controlled frame associated to \,$\left(\,a_{\,2},\, \cdots,\, a_{\,n}\,\right)$\, under some sufficient conditions.

\section{Preliminaries}

\begin{theorem}\cite{Christensen}\label{th1}
Let \,$H_{\,1},\, H_{\,2}$\; be two Hilbert spaces and \;$U \,:\, H_{\,1} \,\to\, H_{\,2}$\; be a bounded linear operator with closed range \;$\mathcal{R}_{\,U}$.\;Then there exists a bounded linear operator \,$U^{\dagger} \,:\, H_{\,2} \,\to\, H_{\,1}$\, such that \,$U\,U^{\dagger}\,x \,=\, x\; \;\forall\; x \,\in\, \mathcal{R}_{\,U}$.
\end{theorem}

The operator \,$U^{\dagger}$\, defined in Theorem (\ref{th1}), is called the pseudo-inverse of \,$U$.

\begin{theorem}\cite{Kreyzig}\label{th1.051}
The set \,$\mathcal{S}\,(\,H_{1}\,)$\; of all self-adjoint operators on a Hilbert space \,$H_{1}$\; is a partially ordered set with respect to the partial order \,$\leq$\, which is defined as for \,$T,\,S \,\in\, \mathcal{S}\,(\,H_{1}\,)$ 
\[T \,\leq\, S \,\Leftrightarrow\, \left<\,T\,f,\, f\,\right> \,\leq\, \left<\,S\,f,\, f\,\right>\; \;\forall\; f \,\in\, H_{1}.\] 
\end{theorem}

\begin{definition}\cite{Kreyzig}
A self-adjoint operator \,$U \,:\, H_{1} \,\to\, H_{1}$\, is called positive if \,$\left<\,U\,x \,,\,  x\,\right> \,\geq\, 0$\, for all \,$x \,\in\, H_{1}$.\;In notation, we can write \,$U \,\geq\, 0$.\;A self-adjoint operator \,$V \,:\, H_{1} \,\to\, H_{1}$\, is called a square root of \,$U$\, if \,$V^{\,2} \,=\, U$.\;If, in addition \,$V \,\geq\, 0$, then \,$V$\, is called positive square root of \,$U$\, and is denoted by \,$V \,=\, U^{1 \,/\, 2}$. 
\end{definition}

\begin{theorem}\cite{Kreyzig}\label{th1.05}
The positive square root \,$V \,:\, H_{1} \,\to\, H_{1}$\, of an arbitrary positive self-adjoint operator \,$U \,:\, H_{1} \,\to\, H_{1}$\, exists and is unique.\;Further, the operator \,$V$\, commutes with every bounded linear operator on \,$H_{1}$\, which commutes with \,$U$.
\end{theorem}

In a complex Hilbert space, every bounded positiove operator is self-adjoint and any two bounded positive operators can be commute with each other.

\begin{definition}\cite{Christensen}
A sequence \,$\left\{\,f_{\,i}\,\right\}_{i \,=\, 1}^{\infty}$\, in a separable Hilbert space \,$H_{1}$\, is said to be a frame for \,$H_{1}$\, if there exist positive constants \,$A,\, B$\, such that
\begin{equation}\label{ee1}
A\; \|\,f\,\|^{\,2} \,\leq\, \sum\limits_{i \,=\, 1}^{\infty}\, \left|\ \left <\,f,\, f_{\,i} \, \right >\,\right|^{\,2} \,\leq\, B \,\|\,f\,\|^{\,2}\; \;\forall\; f \,\in\, H_{1}.
\end{equation}
The constants \,$A$\, and \,$B$\, are called frame bounds.\,If the collection \,$\left\{\,f_{\,i}\,\right\}_{i \,=\, 1}^{\infty}$\, satisfies only the right inequality of (\ref{ee1}) then it is called a Bessel sequence with bound \,$B$.
\end{definition}

\begin{definition}\cite{B}
Let \,$C$\, be a bounded linear operator on \,$H_{1}$\, which has bounded inverse.\,A frame controlled by the operator \,$C$\, or \,$C$-controlled frame is a family of vectors \,$\left\{\,f_{\,i}\,\right\}_{i \,=\, 1}^{\infty}$\, in \,$H_{1}$, such that there exist constants \,$0 \,<\, A \,\leq\, B \,<\, \infty$, satisfying
\[A\; \|\,f\,\|^{\,2} \,\leq\, \sum\limits_{i \,=\, 1}^{\infty}\, \left <\,f,\, f_{\,i} \, \right >\,\left<\,C\,f_{\,i},\, f\,\right> \,\leq\, B \,\|\,f\,\|^{\,2}\; \;\forall\; f \,\in\, H_{1}.\] 
\end{definition}

The controlled frame operator \,$S \,:\, H_{1} \,\to\, H_{1}$\, is defined by
\[S\,f \,=\, \sum\limits_{i \,=\, 1}^{\infty}\, \left <\,f,\, f_{\,i} \, \right >\,C\,f_{\,i}\; \;\forall\; f \,\in\, H_{1}.\]

\begin{definition}\cite{Upender}\label{def0.001}
The tensor product of Hilbert spaces \,$\left(\,H_{1},\, \left<\,,\, \cdot,\,\right>_{1}\,\right)$\, and \,$\left(\,H_{2},\, \left<\,,\, \cdot,\,\right>_{2}\,\right)$\,  is denoted by \,$H_{1} \,\otimes\, H_{2}$\, and it is defined to be an inner product space associated with the inner product   
\[\left<\,f \,\otimes\, g,\, f^{\,\prime} \,\otimes\, g^{\,\prime}\,\right> \,=\, \left<\,f,\, f^{\,\prime}\,\right>_{1}\;\left<\,g,\, g^{\,\prime}\,\right>_{2}\; \;\forall\; f,\, f^{\,\prime} \,\in\, H_{1}\; \;\&\; \;g,\, g^{\,\prime} \,\in\, H_{2}.\]
The norm on \,$H_{1} \,\otimes\, H_{2}$\, is given by 
\[\left\|\,f \,\otimes\, g\,\right\| \,=\, \|\,f\,\|_{\,1}\;\|\,g\,\|_{\,2}\; \;\forall\; f \,\in\, H_{1}\; \;\&\; \,g \,\in\, H_{2}.\]
The space \,$H_{1} \,\otimes\, H_{2}$\, is a Hilbert space with respect to the above inner product.     
\end{definition} 

For \,$Q \,\in\, \mathcal{B}\,(\,H_{1}\,)$\, and \,$T \,\in\, \mathcal{B}\,(\,H_{2}\,)$, the tensor product of operators \,$Q$\, and \,$T$\, is denoted by \,$Q \,\otimes\, T$\, and defined as 
\[\left(\,Q \,\otimes\, T\,\right)\,A \,=\, Q\,A\,T^{\,\ast}\; \;\forall\; \;A \,\in\, H_{1} \,\otimes\, H_{2}.\]

\begin{theorem}\cite{Folland}\label{th1.1}
Suppose \,$Q,\, Q^{\prime} \,\in\, \mathcal{B}\,(\,H_{1}\,)$\, and \,$T,\, T^{\prime} \,\in\, \mathcal{B}\,(\,H_{2}\,)$.\,Then \begin{itemize}
\item[$(i)$]\,$Q \,\otimes\, T \,\in\, \mathcal{B}\,(\,H_{1} \,\otimes\, H_{2}\,)$\, and \,$\left\|\,Q \,\otimes\, T\,\right\| \,=\, \|\,Q\,\|\; \|\,T\,\|$.
\item[$(ii)$] \,$\left(\,Q \,\otimes\, T\,\right)\,(\,f \,\otimes\, g\,) \,=\, Q\,f \,\otimes\, T\,g$\, for all \,$f \,\in\, H_{1},\, g \,\in\, H_{2}$.
\item[$(iii)$] $\left(\,Q \,\otimes\, T\,\right)\,\left(\,Q^{\,\prime} \,\otimes\, T^{\,\prime}\,\right) \,=\, (\,Q\,Q^{\,\prime}\,) \,\otimes\, (\,T\,T^{\,\prime}\,)$. 
\item[$(iv)$] \,$Q \,\otimes\, T$\, is invertible if and only if \,$Q$\, and \,$T$\, are invertible, in which case \,$\left(\,Q \,\otimes\, T\,\right)^{\,-\, 1} \,=\, \left(\,Q^{\,-\, 1} \,\otimes\, T^{\,-\, 1}\,\right)$.
\item[$(v)$] \,$\left(\,Q \,\otimes\, T\,\right)^{\,\ast} \,=\, \left(\,Q^{\,\ast} \,\otimes\, T^{\,\ast}\,\right)$.    
\end{itemize}
\end{theorem}

\begin{definition}\cite{Mashadi}
A \,$n$-norm on a linear space \,$X$\, (\,over the field \,$\mathbb{K}$\, of real or complex numbers\,) is a function
\[\left(\,x_{\,1},\, x_{\,2},\, \cdots,\, x_{\,n}\,\right) \,\longmapsto\, \left\|\,x_{\,1},\, x_{\,2},\, \cdots,\, x_{\,n}\,\right\|,\; x_{\,1},\, x_{\,2},\, \cdots,\, x_{\,n} \,\in\, X\]from \,$X^{\,n}$\, to the set \,$\mathbb{R}$\, of all real numbers such that for every \,$x_{\,1},\, x_{\,2},\, \cdots,\, x_{\,n} \,\in\, X$\, and \,$\alpha \,\in\, \mathbb{K}$,
\begin{itemize}
\item[(i)]\;\; $\left\|\,x_{\,1},\, x_{\,2},\, \cdots,\, x_{\,n}\,\right\| \,=\, 0$\; if and only if \,$x_{\,1},\, \cdots,\, x_{\,n}$\; are linearly dependent,
\item[(ii)]\;\;\; $\left\|\,x_{\,1},\, x_{\,2},\, \cdots,\, x_{\,n}\,\right\|$\; is invariant under permutations of \,$x_{\,1},\, x_{\,2},\, \cdots,\, x_{\,n}$,
\item[(iii)]\;\;\; $\left\|\,\alpha\,x_{\,1},\, x_{\,2},\, \cdots,\, x_{\,n}\,\right\| \,=\, |\,\alpha\,|\, \left\|\,x_{\,1},\, x_{\,2},\, \cdots,\, x_{\,n}\,\right\|$,
\item[(iv)]\;\; $\left\|\,x \,+\, y,\, x_{\,2},\, \cdots,\, x_{\,n}\,\right\| \,\leq\, \left\|\,x,\, x_{\,2},\, \cdots,\, x_{\,n}\,\right\| \,+\,  \left\|\,y,\, x_{\,2},\, \cdots,\, x_{\,n}\,\right\|$.
\end{itemize}
A linear space \,$X$, together with a n-norm \,$\left\|\,\cdot,\, \cdots,\, \cdot \,\right\|$, is called a linear n-normed space. 
\end{definition}

\begin{definition}\cite{Misiak}
Let \,$n \,\in\, \mathbb{N}$\; and \,$X$\, be a linear space of dimension greater than or equal to \,$n$\; over the field \,$\mathbb{K}$, where \,$\mathbb{K}$\, is the real or complex numbers field.\;A n-inner product on \,$X$\, is a map 
\[\left(\,x,\, y,\, x_{\,2},\, \cdots,\, x_{\,n}\,\right) \,\longmapsto\, \left<\,x,\, y \,|\, x_{\,2},\, \cdots,\, x_{\,n} \,\right>,\; x,\, y,\, x_{\,2},\, \cdots,\, x_{\,n} \,\in\, X\]from \,$X^{n \,+\, 1}$\, to the set \,$\mathbb{K}$\, such that for every \,$x,\, y,\, x_{\,1},\, x_{\,2},\, \cdots,\, x_{\,n} \,\in\, X$\, and \,$\alpha \,\in\, \mathbb{K}$,
\begin{itemize}
\item[(i)]\;\; $\left<\,x_{\,1},\, x_{\,1} \,|\, x_{\,2},\, \cdots,\, x_{\,n} \,\right> \,\geq\,  0$\; and \;$\left<\,x_{\,1},\, x_{\,1} \;|\; x_{\,2},\, \cdots,\, x_{\,n} \,\right> \;=\;  0$\; if and only if \;$x_{\,1},\, x_{\,2},\, \cdots,\, x_{\,n}$\; are linearly dependent,
\item[(ii)]\;\; $\left<\,x,\, y \;|\; x_{\,2},\, \cdots,\, x_{\,n} \,\right> \;=\; \left<\,x,\, y \;|\; x_{\,i_{\,2}},\, \cdots,\, x_{\,i_{\,n}} \,\right> $\; for every permutations \\$\left(\, i_{\,2},\, \cdots,\, i_{\,n} \,\right)$\; of \;$\left(\, 2,\, \cdots,\, n \,\right)$,
\item[(iii)]\;\; $\left<\,x,\, y \;|\; x_{\,2},\, \cdots,\, x_{\,n} \,\right> \;=\; \overline{\left<\,y,\, x \;|\; x_{\,2},\, \cdots,\, x_{\,n} \,\right> }$,
\item[(iv)]\;\; $\left<\,\alpha\,x,\, y \;|\; x_{\,2},\, \cdots,\, x_{\,n} \,\right> \;=\; \alpha \,\left<\,x,\, y \;|\; x_{\,2},\, \cdots,\, x_{\,n} \,\right> $,
\item[(v)]\;\; $\left<\,x \,+\, y,\, z \;|\; x_{\,2},\, \cdots,\, x_{\,n} \,\right> \;=\; \left<\,x,\, z \;|\; x_{\,2},\, \cdots,\, x_{\,n} \,\right> \,+\,  \left<\,y,\, z \;|\; x_{\,2},\, \cdots,\, x_{\,n} \,\right>$.
\end{itemize}
A linear space \,$X$\, together with n-inner product \,$\left<\,\cdot,\, \cdot \,|\, \cdot,\, \cdots,\, \cdot\,\right>$\, is called n-inner product space.
\end{definition}

\begin{theorem}\cite{Gunawan}
For \,$n$-inner product space \,$\left(\,X,\, \left<\,\cdot,\, \cdot \,|\, \cdot,\, \cdots,\, \cdot\,\right>\,\right)$, 
\[\left|\,\left<\,x,\, y \,|\, x_{\,2},\,  \cdots,\, x_{\,n}\,\right>\,\right| \,\leq\, \left\|\,x,\, x_{\,2},\, \cdots,\, x_{\,n}\,\right\|\, \left\|\,y,\, x_{\,2},\, \cdots,\, x_{\,n}\,\right\|\]
hold for all \,$x,\, y,\, x_{\,2},\, \cdots,\, x_{\,n} \,\in\, X$, where \[\left \|\,x_{\,1},\, x_{\,2},\, \cdots,\, x_{\,n}\,\right\| \,=\, \sqrt{\left <\,x_{\,1},\, x_{\,1} \;|\; x_{\,2},\,  \cdots,\, x_{\,n}\,\right>},\] is called Cauchy-Schwarz inequality.
\end{theorem}

\begin{definition}\cite{Mashadi}
Let \,$\left(\,X,\, \left\|\,\cdot,\, \cdots,\, \cdot \,\right\|\,\right)$\; be a linear n-normed space.\;A sequence \,$\{\,x_{\,k}\,\}$\; in \,$X$\, is said to converge to some \,$x \,\in\, X$\; if 
\[\lim\limits_{k \to \infty}\,\left\|\,x_{\,k} \,-\, x,\, e_{\,2},\, \cdots,\, e_{\,n} \,\right\| \,=\, 0\]
for every \,$ e_{\,2},\, \cdots,\, e_{\,n} \,\in\, X$\, and it is called a Cauchy sequence if 
\[\lim\limits_{l,\, k \,\to\, \infty}\,\left \|\,x_{l} \,-\, x_{\,k},\, e_{\,2},\, \cdots,\, e_{\,n}\,\right\| \,=\, 0\]
for every \,$ e_{\,2},\, \cdots,\, e_{\,n} \,\in\, X$.\;The space \,$X$\, is said to be complete if every Cauchy sequence in this space is convergent in \,$X$.\;A n-inner product space is called n-Hilbert space if it is complete with respect to its induce norm.
\end{definition}

\begin{definition}\label{def0.1}\cite{Prasenjit}
A sequence \,$\left\{\,f_{\,i}\,\right\}^{\,\infty}_{\,i \,=\, 1}$\, in a \,$n$-Hilbert space \,$H$\, is said to be a frame associated to \,$\left(\,a_{\,2},\, \cdots,\, a_{\,n}\,\right)$\, for \,$H$\, if there exist constant \,$0 \,<\, A \,\leq\, B \,<\, \infty$\, such that
\begin{equation}\label{eee1}
A \, \left\|\,f,\, a_{\,2},\, \cdots,\, a_{\,n} \,\right\|^{\,2} \,\leq\, \sum\limits^{\infty}_{i \,=\, 1}\,\left|\,\left<\,f,\, f_{\,i} \,|\, a_{\,2},\, \cdots,\, a_{\,n}\,\right>\,\right|^{\,2} \,\leq\, B\, \left\|\,f,\, a_{\,2},\, \cdots,\, a_{\,n}\,\right\|^{\,2}
\end{equation} 
for all \,$f \,\in\, H$.\;The infimum of all such \,$B$\, is called the optimal upper frame bound and supremum of all such \,$A$\, is called the optimal lower frame bound.\;A sequence \,$\left\{\,f_{\,i}\,\right\}^{\,\infty}_{\,i \,=\, 1}$\, satisfies only the right inequality of (\ref{eee1}) is called a Bessel sequence associated to \,$\left(\,a_{\,2},\, \cdots,\, a_{\,n}\,\right)$\, in \,$H$\, with bound \,$B$.
\end{definition}

Let \,$L_{F}$\, denote the linear subspace of a \,$n$-Hilbert space \,$H$\, spanned by the non-empty finite set \,$F \,=\, \left\{\,\,a_{\,2},\, a_{\,3},\, \cdots,\, a_{\,n}\,\right\}$, where \,$a_{\,2},\, a_{\,3},\, \cdots,\, a_{\,n}$\, are fixed elements in \,$H$.\;Then the quotient space \,$H \,/\, L_{F}$\, is a normed linear space with respect to the norm, 
\[\left\|\,x \,+\, L_{F}\,\right\|_{F} \,=\, \left\|\,x \,,\, a_{\,2} \,,\,  \cdots \,,\, a_{\,n}\,\right\|\; \;\text{for every}\; x \,\in\, H.\]
Let \,$M_{F}$\, be the algebraic complement of \,$L_{F}$, then \,$H \,=\, L_{F} \,\oplus\, M_{F}$.\;Define   
\[\left<\,x \,,\, y\,\right>_{F} \,=\, \left<\,x \,,\, y \;|\; a_{\,2} \,,\,  \cdots \,,\, a_{\,n}\,\right>\; \;\text{on}\; \;H.\]
Then \,$\left<\,\cdot \,,\, \cdot\,\right>_{F}$\, is a semi-inner product on \,$H$\, and this semi-inner product induces an inner product on the quotient space \,$H \,/\, L_{F}$\; which is given by
\[\left<\,x \,+\, L_{F} \,,\, y \,+\, L_{F}\,\right>_{F} \,=\, \left<\,x \,,\, y\,\right>_{F} \,=\, \left<\,x \,,\, y \,|\, a_{\,2} \,,\,  \cdots \,,\, a_{\,n} \,\right>\;\; \;\forall \;\; x,\, y \,\in\, H.\]
By identifying \,$H \,/\, L_{F}$\; with \,$M_{F}$\; in an obvious way, we obtain an inner product on \,$M_{F}$.\;Now, for every \,$x \,\in\, M_{F}$, we define \,$\|\,x\,\|_{F} \;=\; \sqrt{\left<\,x \,,\, x \,\right>_{F}}$\, and it can be easily verify that \,$\left(\,M_{F} \,,\, \|\,\cdot\,\|_{F}\,\right)$\; is a norm space.\;Let \,$H_{F}$\; be the completion of the inner product space \,$M_{F}$.\\

\begin{definition}\label{def0.0001}\cite{Prasenjit}
Let \,$\left\{\,f_{\,i}\,\right\}_{i \,=\, 1}^{\infty}$\; be a frame associated to \,$\left(\,a_{\,2},\, \cdots,\, a_{\,n}\,\right)$\, for \,$H$.\;Then the bounded linear operator \,$ T_{F} \,:\, l^{\,2}\,(\,\mathbb{N}\,) \,\to\, H_{F}$, defined by \,$T_{F}\,\{\,c_{i}\,\} \,=\, \sum\limits_{i \,=\, 1}^{\infty} \;c_{\,i}\,f_{\,i}$, is called  pre-frame operator and its adjoint operator described by
\[T_{F}^{\,\ast} \,:\, H_{F} \,\to\, l^{\,2}\,(\,\mathbb{N}\,),\;T_{F}^{\,\ast}\,f \,=\, \left \{\, \left <\,f,\, f_{i} \,|\, a_{\,2},\, \cdots,\, a_{\,n}\,\right >\,\right \}_{i \,=\, 1}^{\infty}\] is called the analysis operator.\;The operator \,$S_{F} \,:\, H_{F} \,\to\, H_{F}$\, given by 
\[S_{F}\,f \,=\, T_{F}\,T_{F}^{\,\ast}\,f \,=\, \sum\limits^{\infty}_{i \,=\, 1}\; \left <\,f,\, f_{\,i} \,|\, a_{\,2},\, \cdots,\, a_{\,n} \, \right >\,f_{\,i},\; \;\text{for all}\; \,f \,\in\, H_{F},\] is called the frame operator.
\end{definition}

It is easy to verify that \,$A\,I_{F} \,\leq\, S_{F} \,\leq\, B\,I_{F}$.\,Since \,$S^{\,-1}_{F}$\; commutes with both \,$S_{F}$\; and \,$I_{F}$, multiplying in the inequality, \,$ A\,I_{F} \,\leq\, S_{F} \,\leq\, B\,I_{F}$\, by \,$S^{\,-1}_{F}$, we get \,$B^{\,-1}\,I_{F} \,\leq\, S^{\,-1}_{F} \,\leq\, A^{\,-1}\,I_{F}$.\,For more details on frames in \,$n$-Hilbert spaces and their tesnor product one can go through the papers \cite{Prasenjit, GP, PK}.\\ 

For the remaining part of this paper, \,$\left(\,H,\, \left<\,\cdot,\, \cdot \,|\, \cdot,\, \cdots,\, \cdot\,\right> \,\right)$\; is consider to be a \,$n$-Hilbert space.\,$I_{F}$\, will denote the identity operator on \,$H_{F}$\, and \,$\mathcal{B}\,(\,H_{F}\,)$\, denote the space of all bounded linear operator on \,$H_{F}$.\,$\mathcal{G}\,\mathcal{B}\,(\,H_{F}\,)$\, denotes the set of all bounded linear operators which have bounded inverse.\,If \,$S,\, R \,\in\, \mathcal{G}\,\mathcal{B}\,(\,H_{F}\,)$, then \,$R^{\,\ast},\, R^{\,-\, 1}$\, and \,$S\,R$\, are also belongs to \,$\mathcal{G}\,\mathcal{B}\,(\,H_{F}\,)$.\,$\mathcal{G}\,\mathcal{B}^{\,+}\,(\,H_{F}\,)$\, is the set of all positive operators in \,$\mathcal{G}\,\mathcal{B}\,(\,H_{F}\,)$\, and \,$C,\,C_{1},\,C_{2}$\, are invertible operators in \,$\mathcal{G}\,\mathcal{B}\,\left(\,H_{F}\,\right)$.\,$l^{\,2}(\,\mathbb{N}\,)$\, denote the space of square summable scalar-valued sequences with index set of natural numbers \,$\mathbb{N}$.

\section{Controlled frame in $n$-Hilbert space}

\smallskip\hspace{.6 cm} In this section, we introduce the notion of controlled frame in \,$n$-Hilbert space and discuss it's several properties.\,At the end, dual controlled frame in \,$n$-Hilbert space is presented.

\begin{definition}
Let \,$C \,\in\, \mathcal{G}\,\mathcal{B}\left(\,H_{F}\,\right)$.\,A frame associated to \,$\left(\,a_{\,2},\, \cdots,\, a_{\,n}\,\right)$\, controlled by the operator \,$C$\, or \,$C$-controlled frame associated to \,$\left(\,a_{\,2},\, \cdots,\, a_{\,n}\,\right)$\, is a family of vectors \,$\left\{\,f_{\,i}\,\right\}^{\,\infty}_{\,i \,=\, 1}$\, in \,$H$\, such that there exist constants \,$0 \,<\, A \,\leq\, B \,<\, \infty$\, satisfying
\begin{align}
A \, \left\|\,f,\, a_{\,2},\, \cdots,\, a_{\,n} \,\right\|^{\,2}& \,\leq\, \sum\limits_{i \,=\, 1}^{\,\infty}\,\left<\,f,\, f_{\,i} \,|\, a_{\,2},\, \cdots,\, a_{\,n}\,\right>\,\left<\,C\,f_{\,i},\, f \,|\, a_{\,2},\, \cdots,\, a_{\,n}\,\right>\nonumber\\
&\,\leq\, B\, \left\|\,f,\, a_{\,2},\, \cdots,\, a_{\,n}\,\right\|^{\,2}\; \;\forall\; f \,\in\, H_{F}\label{eq1}.
\end{align}  
\end{definition}
The constants \,$A$\, and \,$B$\, are called lower and upper bounds of \,$C$-controlled frame associated to \,$\left(\,a_{\,2},\, \cdots,\, a_{\,n}\,\right)$, respectively. 

\begin{description}
\item[$(i)$]The family of vectors \,$\left\{\,f_{\,i}\,\right\}^{\,\infty}_{\,i \,=\, 1}$\, is called a \,$C$-controlled tight frame associated to \,$\left(\,a_{\,2},\, \cdots,\, a_{\,n}\,\right)$\, if \,$A \,=\, B$\,  and it is called \,$C$-controlled Parseval frame associated to \,$\left(\,a_{\,2},\, \cdots,\, a_{\,n}\,\right)$\, if \,$A \,=\, B \,=\, 1$.
\item[$(ii)$]If only the right inequality of (\ref{eq1}) is satisfied then it is called a \,$C$-controlled Bessel sequence associated to \,$\left(\,a_{\,2},\, \cdots,\, a_{\,n}\,\right)$\, in \,$H$\, with bound \,$B$.  
\end{description}

\begin{remark}
Suppose that \,$\left\{\,f_{\,i}\,\right\}^{\,\infty}_{\,i \,=\, 1}$\, is a \,$C$-controlled tight frame associated to \,$\left(\,a_{\,2},\, \cdots,\, a_{\,n}\,\right)$\, for \,$H$\, with bounds \,$A$.\,Then for each \,$f \,\in\, H_{F}$, we have
\begin{align*}
&\sum\limits_{i \,=\, 1}^{\,\infty}\,\left<\,f,\, f_{\,i} \,|\, a_{\,2},\, \cdots,\, a_{\,n}\,\right>\,\left<\,C\,f_{\,i},\, f \,|\, a_{\,2},\, \cdots,\, a_{\,n}\,\right> \,=\, A\, \left\|\,f,\, a_{\,2},\, \cdots,\, a_{\,n}\,\right\|^{\,2}\\
&\Rightarrow\,\sum\limits_{i \,=\, 1}^{\,\infty}\,\left<\,f,\, A^{\,-\, 1 \,/\, 2}\,f_{\,i} \,|\, a_{\,2},\, \cdots,\, a_{\,n}\,\right>\,\left<\,C\,A^{\,-\, 1 \,/\, 2}\,f_{\,i},\, f \,|\, a_{\,2},\, \cdots,\, a_{\,n}\,\right>\\
&\hspace{1.3cm} \,=\, \left\|\,f,\, a_{\,2},\, \cdots,\, a_{\,n}\,\right\|^{\,2}. 
\end{align*} 
This verify that \,$\left\{\,A^{\,-\, 1 \,/\, 2}\,f_{\,i}\,\right\}^{\,\infty}_{\,i \,=\, 1}$\, is a \,$C$-controlled Parseval frame associated to \,$\left(\,a_{\,2},\, \cdots,\, a_{\,n}\,\right)$\, for \,$H$. 
\end{remark}

\begin{definition}\label{def0.0001}
Let \,$\left\{\,f_{\,i}\,\right\}^{\,\infty}_{\,i \,=\, 1}$\, be a \,$C$-controlled Bessel sequence associated to \,$\left(\,a_{\,2},\, \cdots,\, a_{\,n}\,\right)$\, in \,$H$.\,The operator \,$T_{C} \,:\, l^{\,2}\,(\,\mathbb{N}\,) \,\to\, H_{F}$\, given by \,$T_{C}\,\left(\,\left\{\,c_{i}\,\right\}_{i \,=\, 1}^{\,\infty}\,\right) \,=\, \sum\limits^{\infty}_{i \,=\, 1}\;c_{\,i}\,C\,f_{\,i}$\, is called pre-frame operator or the synthesis operator.\,The adjoint operator described by 
\[T^{\,\ast}_{C} \,:\, H_{F} \,\to\, l^{\,2}\,(\,\mathbb{N}\,),\; T_{C}^{\,\ast}\,f \,=\, \left\{\,\left<\,f,\, f_{\,i} \,|\, a_{\,2},\, \cdots,\, a_{\,n} \,\right>\,\right\}_{i \,\in\, I},\, \,f \,\in\, H_{F}\]
is called the analysis operator.\,The \,$C$-controlled frame operator \,$S_{C} \,:\, H_{F} \,\to\, H_{F}$\, is defined as 
\[S_{C}\,f \,=\, T_{C}\,T^{\,\ast}_{C}\,f \,=\, \sum\limits_{i \,=\, 1}^{\,\infty}\, \left<\,f,\, f_{\,i} \,|\, a_{\,2},\, \cdots,\, a_{\,n} \,\right>\,C\,f_{\,i}\; \;\forall\; f \,\in\, H_{F}.\]
\end{definition}

For each \,$f \,\in\, H_{F}$, we have
\begin{align*}
\left<\,S_{C}\,f,\, f \,|\, a_{\,2},\, \cdots,\, a_{\,n}\,\right> &\,=\, \left<\,\sum\limits_{i \,=\, 1}^{\,\infty}\, \left<\,f,\, f_{\,i} \,|\, a_{\,2},\, \cdots,\, a_{\,n} \,\right>\,C\,f_{\,i},\, f \,|\, a_{\,2},\, \cdots,\, a_{\,n}\,\right>\\
&=\, \sum\limits_{i \,=\, 1}^{\,\infty}\, \left<\,f,\, f_{\,i} \,|\, a_{\,2},\, \cdots,\, a_{\,n} \,\right>\,\left<\,C\,f_{\,i},\, f \,|\, a_{\,2},\, \cdots,\, a_{\,n} \,\right>. 
\end{align*}
Thus, if \,$\left\{\,f_{\,i}\,\right\}^{\,\infty}_{\,i \,=\, 1}$\, is a \,$C$-controlled frame associated to \,$\left(\,a_{\,2},\, \cdots,\, a_{\,n}\,\right)$\, for \,$H$\, then from (\ref{eq1}), we have
\begin{align*}
&A\,\left<\,f,\, f \,|\, a_{\,2},\, \cdots,\, a_{\,n}\,\right> \,\leq\, \left<\,S_{C}\,f,\, f \,|\, a_{\,2},\, \cdots,\, a_{\,n}\,\right> \,\leq\, B\,\left<\,f,\, f \,|\, a_{\,2},\, \cdots,\, a_{\,n}\,\right>\\
&\Rightarrow\, A\,I_{H} \,\leq\, S_{C} \,\leq\, B\,I_{H}.
\end{align*}
From this we can conclude that \,$S_{C}$\, is bounded, self-adjoint, positive and invertible.

\begin{theorem}\label{th2.1}
Let \,$\left\{\,f_{\,i}\,\right\}^{\,\infty}_{\,i \,=\, 1}$\, be a \,$C$-controlled Bessel sequence associated to \,$\left(\,a_{\,2},\, \cdots,\, a_{\,n}\,\right)$\, in \,$H$\, with bound \,$B$.\,Then the operator given by
\[U \,:\, l^{\,2}\,(\,\mathbb{N}\,) \,\to\, H_{F},\; U\,\left(\,\left\{\,a_{i}\,\right\}_{i \,=\, 1}^{\,\infty}\,\right) \,=\, \sum\limits^{\infty}_{i \,=\, 1}\;a_{\,i}\,C\,f_{\,i}\] is well-defined and bounded operator from \,$l^{\,2}\,(\,\mathbb{N}\,)$\, into \,$H_{F}$\, with \,$\|\,U\,\| \,\leq\, \sqrt{B}\,\left\|\,C^{\,1 \,/\, 2}\,\right\|$.
\end{theorem}

\begin{proof}
Suppose \,$\left\{\,f_{\,i}\,\right\}^{\,\infty}_{\,i \,=\, 1}$\, is a \,$C$-controlled Bessel sequence associated to \,$\left(\,a_{\,2},\, \cdots,\, a_{\,n}\,\right)$\, in \,$H$\, with bound \,$B$.\,Then for each \,$f \,\in\, H_{F}$, we have
\[\sum\limits_{i \,=\, 1}^{\,\infty}\,\left<\,f,\, f_{\,i} \,|\, a_{\,2},\, \cdots,\, a_{\,n}\,\right>\,\left<\,C\,f_{\,i},\, f \,|\, a_{\,2},\, \cdots,\, a_{\,n}\,\right> \,\leq\, B\,\left\|\,f,\, a_{\,2},\, \cdots,\, a_{\,n} \,\right\|^{\,2}.\]
Let \,$\left\{\,a_{i}\,\right\}_{i \,=\, 1}^{\,\infty} \,\in\, l^{\,2}\,(\,\mathbb{N}\,)$.\,For arbitrary \,$l \,>\, k$, we have
\begin{align*}
&\left\|\,\sum\limits^{l}_{i \,=\, 1}\,a_{\,i}\,C\,f_{\,i} \,-\, \sum\limits^{k}_{i \,=\, 1}\,a_{i}\,C\,f_{\,i}\,\right\|_{F}^{\,2} \,=\, \left\|\,\sum\limits^{l}_{i \,=\, k+1}\,a_{\,i}\,C\,f_{\,i}\, \;,\; a_{\,2} \,,\, \cdots \,,\, a_{\,n}\,\right\|^{\,2}\\
&=\sup\limits_{f \,\in\, H_{F},\,\left\|\,f,\, a_{\,2},\, \cdots,\, a_{\,n}\,\right\| \,=\, 1}\,\,\left\{\,\left|\,\left<\,\sum\limits^{l}_{i \,=\, k+1}\,a_{\,i}\,C\,f_{\,i}, f \,|\, a_{\,2},\, \cdots,\, a_{\,n}\,\right>\,\right|^{\,2}\,\right\}\\
&=\,\sup\limits_{f \,\in\, H_{F},\,\left\|\,f,\, a_{\,2},\, \cdots,\, a_{\,n}\,\right\| \,=\, 1}\,\left\{\,\left|\,\sum\limits^{l}_{i \,=\, k+1}\,a_{\,i}\,\left<\,C\,f_{\,i}, f \,|\, a_{\,2},\, \cdots,\, a_{\,n}\,\right>\,\right|^{\,2}\,\right\}\\
&\leq\,\sup\limits_{f \,\in\, H_{F},\,\left\|\,f,\, a_{\,2},\, \cdots,\, a_{\,n}\,\right\| \,=\, 1}\,\left\{\,\sum\limits^{l}_{i \,=\, k+1}\,\left<\,f,\, C\,f_{\,i} \,|\, a_{\,2},\, \cdots,\, a_{\,n}\,\right>\,\left<\,C\,f_{\,i},\, f \,|\, a_{\,2},\, \cdots,\, a_{\,n}\,\right>\,\right\}\,\times\\
&\hspace{2cm}\sum\limits_{i \,=\, k \,+\, 1}^{\,l}\,\left|\,a_{i}\,\right|^{\,2}\\
&\leq\sup\limits_{\left\|\,f,\, a_{\,2},\, \cdots,\, a_{\,n}\,\right\| \,=\, 1}\left<\,\sum\limits_{i \,=\, k \,+\, 1}^{\,l}\left<\,f,\, C\,f_{\,i} \,|\, a_{\,2},\, \cdots,\, a_{\,n}\,\right>\,C\,f_{\,i},\, f \,|\, a_{\,2},\, \cdots,\, a_{\,n}\,\right>\left\|\,\left\{\,a_{i}\,\right\}_{i \,\in\, I}\,\right\|^{\,2}\\
&=\,\sup\limits_{\left\|\,f,\, a_{\,2},\, \cdots,\, a_{\,n}\,\right\| \,=\, 1}\,\left<\,C\,S_{C}\,f,\, f \,|\, a_{\,2},\, \cdots,\, a_{\,n}\,\right>\,\left\|\,\left\{\,a_{i}\,\right\}_{i \,\in\, I}\,\right\|^{\,2}\\
&=\,\sup\limits_{\left\|\,f,\, a_{\,2},\, \cdots,\, a_{\,n}\,\right\| \,=\, 1}\,\left<\,\left(\,C\,S_{C}\,\right)^{1 \,/\, 2}\,f,\, \left(\,C\,S_{C}\,\right)^{1 \,/\, 2}\,f \,|\, a_{\,2},\, \cdots,\, a_{\,n}\,\right>\,\left\|\,\left\{\,a_{i}\,\right\}_{i \,\in\, I}\,\right\|^{\,2}\\ 
&=\,\sup\limits_{\left\|\,f,\, a_{\,2},\, \cdots,\, a_{\,n}\,\right\| \,=\, 1}\,\left\|\,\left(\,C\,S_{C}\,\right)^{1 \,/\, 2}\,f,\, a_{\,2},\, \cdots,\, a_{\,n}\,\right\|^{\,2}\,\left\|\,\left\{\,a_{i}\,\right\}_{i \,\in\, I}\,\right\|^{\,2}\\ 
&\leq\,\sup\limits_{\left\|\,f,\, a_{\,2},\, \cdots,\, a_{\,n}\,\right\| \,=\, 1}\,\left\|\,C^{1 \,/\, 2}\,\right\|^{\,2}\,\left\|\,S_{C}^{1 \,/\, 2}\,f,\, a_{\,2},\, \cdots,\, a_{\,n}\,\right\|^{\,2}\,\left\|\,\left\{\,a_{i}\,\right\}_{i \,\in\, I}\,\right\|^{\,2}\\
&=\,\sup\limits_{\left\|\,f,\, a_{\,2},\, \cdots,\, a_{\,n}\,\right\| \,=\, 1}\,\left\|\,C^{1 \,/\, 2}\,\right\|^{\,2}\,\left<\,S_{C}\,f,\, f \,|\, a_{\,2},\, \cdots,\, a_{\,n}\,\right>\,\left\|\,\left\{\,a_{i}\,\right\}_{i \,\in\, I}\,\right\|^{\,2}\\
&\leq\, B\,\left\|\,C^{1 \,/\, 2}\,\right\|^{\,2}\,\left\|\,\left\{\,a_{i}\,\right\}_{i \,\in\, I}\,\right\|^{\,2}. 
\end{align*}
This shows that \,$\sum\limits^{\infty}_{i \,=\, 1}\,a_{\,i}\,C\,f_{\,i}$\, is a Cauchy sequence which is convergent in \,$H_{F}$.\,Thus \,$U$\, is well-defined and bounded with \,$\|\,U\,\| \,\leq\, \sqrt{B}\,\left\|\,C^{\,1 \,/\, 2}\,\right\|$.
\end{proof}

For arbitrary \,$f \,\in\, H_{F}$\, and \,$\left\{\,a_{i}\,\right\}_{i \,=\, 1}^{\,\infty} \,\in\, l^{\,2}\,(\,\mathbb{N}\,)$, we have
\begin{align*}
\left<\,f,\, U\,\left\{\,a_{i}\,\right\} \,|\, a_{\,2},\, \cdots,\, a_{\,n}\,\right>& \,=\, \left<\,f,\, \sum\limits_{i \,=\, 1}^{\,\infty}\,a_{i}\,C\,f_{\,i} \,|\, a_{\,2},\, \cdots,\, a_{\,n}\,\right>\\
&=\,\,\sum\limits_{i \,=\, 1}^{\,\infty}\,\overline{\,a_{i}}\,\left<\,C\,f,\, f_{\,i} \,|\, a_{\,2},\, \cdots,\, a_{\,n}\,\right>.
\end{align*} 
Therefore, 
\[\left<\,f,\, U\,\left\{\,a_{i}\,\right\}\right>_{F} \,=\, \left<\,\left\{\,\left<\,C\,f,\, f_{\,i} \,|\, a_{\,2},\, \cdots,\, a_{\,n}\,\right>\,\right\},\, \left\{\,a_{i}\,\right\}\,\right>\]
and hence
\[U^{\,\ast}\,f \,=\, \left\{\,\left<\,C\,f,\, f_{\,i} \,|\, a_{\,2},\, \cdots,\, a_{\,n}\,\right>\,\right\}_{i \,=\, 1}^{\,\infty}\; \;\forall\; f \,\in\, H_{F}.\]  

The following theorem shows that any controlled frame associated to \,$\left(\,a_{\,2},\, \cdots,\, a_{\,n}\,\right)$\, is a frame associated to \,$\left(\,a_{\,2},\, \cdots,\, a_{\,n}\,\right)$. 

\begin{theorem}
Let \,$\left\{\,f_{\,i}\,\right\}^{\,\infty}_{\,i \,=\, 1}$\, be a \,$C$-controlled frame associated to \,$\left(\,a_{\,2},\, \cdots,\, a_{\,n}\,\right)$\, for \,$H$\, and \,$C \,\in\, \mathcal{G}\,\mathcal{B}^{\,+\,}\left(\,H_{F}\,\right)$.\,Then \,$\left\{\,f_{\,i}\,\right\}^{\,\infty}_{\,i \,=\, 1}$\, is a frame associated to \,$\left(\,a_{\,2},\, \cdots,\, a_{\,n}\,\right)$\, for \,$H$.  
\end{theorem}

\begin{proof}
Suppose that \,$\left\{\,f_{\,i}\,\right\}^{\,\infty}_{\,i \,=\, 1}$\, is a \,$C$-controlled frame associated to \,$\left(\,a_{\,2},\, \cdots,\, a_{\,n}\,\right)$\, for \,$H$\, with bounds \,$A$\, and \,$B$.\,Then for each \,$f \,\in\, H_{F}$, we have
\begin{align*}
&A\,\left\|\,f,\, a_{\,2},\, \cdots,\, a_{\,n} \,\right\|^{\,2} \,=\, A\,\left\|\,C^{\, 1 \,/\, 2}\,C^{\,-\, 1 \,/\, 2}\,f,\, a_{\,2},\, \cdots,\, a_{\,n} \,\right\|^{\,2}\\
&\leq\, A\,\left\|\,C^{\,1 \,/\, 2}\,\right\|^{\,2}\,\left\|\,C^{\,-\, 1 \,/\, 2}\,f,\, a_{\,2},\, \cdots,\, a_{\,n} \,\right\|^{\,2}\\
&\leq\, \left\|\,C^{\,1 \,/\, 2}\,\right\|^{\,2}\,\sum\limits_{i \,=\, 1}^{\,\infty}\,\left<\,C^{\,-\, 1 \,/\, 2}\,f,\, f_{\,i} \,|\, a_{\,2},\, \cdots,\, a_{\,n}\,\right>\,\left<\,C\,f_{\,i},\, C^{\,-\, 1 \,/\, 2}\,f \,|\, a_{\,2},\, \cdots,\, a_{\,n}\,\right>\\
&=\, \left\|\,C^{\,1 \,/\, 2}\,\right\|^{\,2}\,\sum\limits_{i \,=\, 1}^{\,\infty}\,\left<\,C^{\,-\, 1 \,/\, 2}\,f,\, f_{\,i} \,|\, a_{\,2},\, \cdots,\, a_{\,n}\,\right>\,\left<\,f_{\,i},\, C^{\, 1 \,/\, 2}\,f \,|\, a_{\,2},\, \cdots,\, a_{\,n}\,\right>\\
&=\, \left\|\,C^{\,1 \,/\, 2}\,\right\|^{\,2}\,\left<\,\sum\limits_{i \,=\, 1}^{\,\infty}\,\left<\,C^{\,-\, 1 \,/\, 2}\,f,\, f_{\,i} \,|\, a_{\,2},\, \cdots,\, a_{\,n}\,\right>\,f_{\,i},\, C^{\, 1 \,/\, 2}\,f \,|\, a_{\,2},\, \cdots,\, a_{\,n}\,\right>\\
&=\, \left\|\,C^{\,1 \,/\, 2}\,\right\|^{\,2}\,\left<\,S_{F}\,C^{\,-\, 1 \,/\, 2}\,f,\, C^{\, 1 \,/\, 2}\,f \,|\, a_{\,2},\, \cdots,\, a_{\,n}\,\right>\\
&=\, \left\|\,C^{\,1 \,/\, 2}\,\right\|^{\,2}\,\left<\,S_{F}\,f,\, f \,|\, a_{\,2},\, \cdots,\, a_{\,n}\,\right> \,=\, \left\|\,C^{\,1 \,/\, 2}\,\right\|^{\,2}\,\sum\limits_{i \,=\, 1}^{\,\infty}\,\left|\,\left<\,f,\, f_{\,i} \,|\, a_{\,2},\, \cdots,\, a_{\,n}\,\right>\,\right|^{\,2}\\
&\Rightarrow\,A\,\left\|\,C^{\,1 \,/\, 2}\,\right\|^{\,-\, 2}\,\left\|\,f,\, a_{\,2},\, \cdots,\, a_{\,n} \,\right\|^{\,2} \,\leq\, \sum\limits_{i \,=\, 1}^{\,\infty}\,\left|\,\left<\,f,\, f_{\,i} \,|\, a_{\,2},\, \cdots,\, a_{\,n}\,\right>\,\right|^{\,2}.   
\end{align*} 
On the other hand, for each \,$f \,\in\, H_{F}$, we have
\begin{align*}
&\sum\limits_{i \,=\, 1}^{\,\infty}\,\left|\,\left<\,f,\, f_{\,i} \,|\, a_{\,2},\, \cdots,\, a_{\,n}\,\right>\,\right|^{\,2} \,=\, \left<\,f,\, S_{F}\,f \,|\, a_{\,2},\, \cdots,\, a_{\,n}\,\right>\\
& =\,\left<\,f,\, C^{\,-\, 1}\,C\,S_{F}\,f \,|\, a_{\,2},\, \cdots,\, a_{\,n}\,\right>\\
&=\, \left<\,\left(\,C^{\,-\, 1}\,C\,S_{F}\,\right)^{\,1 \,/\, 2} f,\, \left(\,C^{\,-\, 1}\,C\,S_{F}\,\right)^{\,1 \,/\, 2}\,f \,|\, a_{\,2},\, \cdots,\, a_{\,n}\,\right>\\
&=\, \left\|\,\left(\,C^{\,-\, 1}\,C\,S_{F}\,\right)^{\,1 \,/\, 2}\,f,\, a_{\,2},\, \cdots,\, a_{\,n}\,\right\|^{\,2}\\
&\leq\,\left\|\,C^{\,-\, 1 \,/\, 2}\,\right\|^{\,2}\,\left\|\,\left(\,C\,S_{F}\,\right)^{\,1 \,/\, 2}\,f,\, a_{\,2},\, \cdots,\, a_{\,n}\,\right\|^{\,2}\\
&=\,\left\|\,C^{\,-\, 1 \,/\, 2}\,\right\|^{\,2}\,\left<\,f,\, C\,S_{F}\,f \,|\, a_{\,2},\, \cdots,\, a_{\,n}\,\right>\\ 
&=\,\left\|\,C^{\,-\, 1 \,/\, 2}\,\right\|^{\,2}\,\left<\,f,\, S_{C}\,f \,|\, a_{\,2},\, \cdots,\, a_{\,n}\,\right>\\ 
&=\,\left\|\,C^{\,-\, 1 \,/\, 2}\,\right\|^{\,2}\,\sum\limits_{i \,=\, 1}^{\,\infty}\,\left<\,f,\, f_{\,i} \,|\, a_{\,2},\, \cdots,\, a_{\,n}\,\right>\,\left<\,C\,f_{\,i},\, f \,|\, a_{\,2},\, \cdots,\, a_{\,n}\,\right>\\
&\leq\, B\,\left\|\,C^{\,-\, 1 \,/\, 2}\,\right\|^{\,2}\,\left\|\,f,\, a_{\,2},\, \cdots,\, a_{\,n} \,\right\|^{\,2}.   
\end{align*}
Thus, \,$\left\{\,f_{\,i}\,\right\}^{\,\infty}_{\,i \,=\, 1}$\, is a frame associated to \,$\left(\,a_{\,2},\, \cdots,\, a_{\,n}\,\right)$\, for \,$H$\, with bounds \\$A\,\left\|\,C^{\,1 \,/\, 2}\,\right\|^{\,-\, 2}$\, and \,$B\,\left\|\,C^{\,-\, 1 \,/\, 2}\,\right\|^{\,2}$.\,This completes the proof. 
\end{proof}

Next theorem shows that any frame associated to \,$\left(\,a_{\,2},\, \cdots,\, a_{\,n}\,\right)$\, is a controlled frame associated to \,$\left(\,a_{\,2},\, \cdots,\, a_{\,n}\,\right)$\, under some conditions.

\begin{theorem}
Let \,$C \,\in\, \mathcal{G}\,\mathcal{B}^{\,+\,}\left(\,H_{F}\,\right)$\, be a self-adjoint operator.\,If \,$\left\{\,f_{\,i}\,\right\}^{\,\infty}_{\,i \,=\, 1}$\, is a frame associated to \,$\left(\,a_{\,2},\, \cdots,\, a_{\,n}\,\right)$\, for \,$H$\, then \,$\left\{\,f_{\,i}\,\right\}^{\,\infty}_{\,i \,=\, 1}$\, is a \,$C$-controlled frame associated to \,$\left(\,a_{\,2},\, \cdots,\, a_{\,n}\,\right)$\, for \,$H$. 
\end{theorem}

\begin{proof}
Suppose that \,$\left\{\,f_{\,i}\,\right\}^{\,\infty}_{\,i \,=\, 1}$\, is a frame associated to \,$\left(\,a_{\,2},\, \cdots,\, a_{\,n}\,\right)$\, for \,$H$\, with bounds \,$A$\, and \,$B$.\,Then for each \,$f \,\in\, H_{F}$, we have
\begin{align*}
&A\,\left\|\,f,\, a_{\,2},\, \cdots,\, a_{\,n} \,\right\|^{\,2} \,=\, A\,\left\|\,C^{\,-\, 1 \,/\, 2}\,C^{\,1 \,/\, 2}\,f,\, a_{\,2},\, \cdots,\, a_{\,n} \,\right\|^{\,2}\\
&\leq\, A\,\left\|\,C^{\,-\, 1 \,/\, 2}\,\right\|^{\,2}\,\left\|\,C^{\,1 \,/\, 2}\,f,\, a_{\,2},\, \cdots,\, a_{\,n} \,\right\|^{\,2}\\
&\leq\,\left\|\,C^{\,-\, 1 \,/\, 2}\,\right\|^{\,2}\,\sum\limits_{i \,=\, 1}^{\,\infty}\,\left<\,C^{\,1 \,/\, 2}\,f,\, f_{\,i} \,|\, a_{\,2},\, \cdots,\, a_{\,n}\,\right>\,\left<\,C^{\,1 \,/\, 2}\,f_{\,i},\, f \,|\, a_{\,2},\, \cdots,\, a_{\,n}\,\right>\\
&=\,\left\|\,C^{\,-\, 1 \,/\, 2}\,\right\|^{\,2}\,\left<\,C^{\,1 \,/\, 2}\,f,\, \sum\limits_{i \,=\, 1}^{\,\infty}\,\left<\,f_{\,i},\, C^{\,1 \,/\, 2}\,f \,|\, a_{\,2},\, \cdots,\, a_{\,n}\,\right>\,f_{\,i} \,|\, a_{\,2},\, \cdots,\, a_{\,n}\,\right>\\
&=\,\left\|\,C^{\,-\, 1 \,/\, 2}\,\right\|^{\,2}\,\left<\,C^{\,1 \,/\, 2}\,f,\, C^{\,1 \,/\, 2}\,S_{F}\,f \,|\, a_{\,2},\, \cdots,\, a_{\,n}\,\right>\\ 
&=\,\left\|\,C^{\,-\, 1 \,/\, 2}\,\right\|^{\,2}\,\left<\,f,\, C\,S_{F}\,f \,|\, a_{\,2},\, \cdots,\, a_{\,n}\,\right>\\
&\Rightarrow\, A\,\left\|\,C^{\,-\, 1 \,/\, 2}\,\right\|^{\,-\, 2}\,\left\|\,f,\, a_{\,2},\, \cdots,\, a_{\,n} \,\right\|^{\,2}\\
& \,\leq\, \sum\limits_{i \,=\, 1}^{\,\infty}\,\left<\,f,\, f_{\,i} \,|\, a_{\,2},\, \cdots,\, a_{\,n}\,\right>\,\left<\,C\,f_{\,i},\, f \,|\, a_{\,2},\, \cdots,\, a_{\,n}\,\right>. 
\end{align*}
On the other hand, for each \,$f \,\in\, H_{F}$, we have
\begin{align*}
&\sum\limits_{i \,=\, 1}^{\,\infty}\,\left<\,f,\, f_{\,i} \,|\, a_{\,2},\, \cdots,\, a_{\,n}\,\right>\,\left<\,C\,f_{\,i},\, f \,|\, a_{\,2},\, \cdots,\, a_{\,n}\,\right>\\
&=\,\left<\,f,\, C\,S_{F}\,f \,|\, a_{\,2},\, \cdots,\, a_{\,n}\,\right> \,=\, \left<\,C\,f,\, S_{F}\,f \,|\, a_{\,2},\, \cdots,\, a_{\,n}\,\right>\\
&\leq\,\left\|\,C\,f,\, a_{\,2},\, \cdots,\, a_{\,n} \,\right\|\,\left\|\,S_{F}\,f,\, a_{\,2},\, \cdots,\, a_{\,n} \,\right\|\\
&\leq\,\|\,C\,\|\,\left\|\,f,\, a_{\,2},\, \cdots,\, a_{\,n} \,\right\|\,B\,\left\|\,f,\, a_{\,2},\, \cdots,\, a_{\,n} \,\right\|\\
&\,=\, B\,\|\,C\,\|\,\left\|\,f,\, a_{\,2},\, \cdots,\, a_{\,n} \,\right\|^{\,2}. 
\end{align*}
Thus, \,$\left\{\,f_{\,i}\,\right\}^{\,\infty}_{\,i \,=\, 1}$\, is a \,$C$-controlled frame associated to \,$\left(\,a_{\,2},\, \cdots,\, a_{\,n}\,\right)$\, for \,$H$. \\This completes the proof. 
\end{proof}

For every controlled frame associated to \,$\left(\,a_{\,2},\, \cdots,\, a_{\,n}\,\right)$\, we can get a canonical controlled tight frame associated to \,$\left(\,a_{\,2},\, \cdots,\, a_{\,n}\,\right)$\, with frame bound \,$1$.

\begin{theorem}
Let \,$\left\{\,f_{\,i}\,\right\}^{\,\infty}_{\,i \,=\, 1}$\, be a \,$C$-controlled frame associated to \,$\left(\,a_{\,2},\, \cdots,\, a_{\,n}\,\right)$\, for \,$H$\, with frame operator \,$S_{C}$.\,Suppose \,$S_{C}^{\,-\, 1}$\, commutes with \,$C$.\,Then \,$\left\{\,S_{C}^{\,-\, 1 \,/\, 2}\,f_{\,i}\,\right\}^{\,\infty}_{\,i \,=\, 1}$\, is a \,$C$-controlled Parseval frame associated to \,$\left(\,a_{\,2},\, \cdots,\, a_{\,n}\,\right)$\, for \,$H$.  
\end{theorem}

\begin{proof}
The existence of unique positive square root of \,$S_{C}^{\,-\, 1}$\, follows from Theorem \ref{th1.051}, since \,$S_{C}^{\,-\, 1 \,/\, 2}$\, is a limit of a sequence of polynomials in \,$S_{C}^{\,-\, 1}$, it commutes with \,$S_{C}^{\,-\, 1}$\, and therefore with \,$S_{C}$.\,Then for each \,$f \,\in\, H_{F}$, we have
\begin{align*}
f &\,=\, S_{C}^{\,-\, 1 \,/\, 2}\,S_{C}\,S_{C}^{\,-\, 1 \,/\, 2}\,f\\
&=\,\sum\limits_{i \,=\, 1}^{\,\infty}\, \left<\,f,\, S_{C}^{\,-\, 1 \,/\, 2}\,f_{\,i} \,|\, a_{\,2},\, \cdots,\, a_{\,n} \,\right>\,S_{C}^{\,-\, 1 \,/\, 2}\,C\,f_{\,i}\\
&=\,\sum\limits_{i \,=\, 1}^{\,\infty}\, \left<\,f,\, S_{C}^{\,-\, 1 \,/\, 2}\,f_{\,i} \,|\, a_{\,2},\, \cdots,\, a_{\,n} \,\right>\,C\,S_{C}^{\,-\, 1 \,/\, 2}\,f_{\,i}.  
\end{align*}
Now, for each \,$f \,\in\, H_{F}$, we have
\begin{align*}
&\left\|\,f,\, a_{\,2},\, \cdots,\, a_{\,n}\,\right\|^{\,2} \,=\, \left<\,f,\, f \,|\, a_{\,2},\, \cdots,\, a_{\,n}\,\right>\\
&\,=\, \left<\,\sum\limits_{i \,=\, 1}^{\,\infty}\, \left<\,f,\, S_{C}^{\,-\, 1 \,/\, 2}\,f_{\,i} \,|\, a_{\,2},\, \cdots,\, a_{\,n} \,\right>\,C\,S_{C}^{\,-\, 1 \,/\, 2}\,f_{\,i},\, f \,|\, a_{\,2},\, \cdots,\, a_{\,n}\,\right>\\
&=\,\sum\limits_{i \,=\, 1}^{\,\infty}\, \left<\,f,\, S_{C}^{\,-\, 1 \,/\, 2}\,f_{\,i} \,|\, a_{\,2},\, \cdots,\, a_{\,n} \,\right>\,\left<\,C\,S_{C}^{\,-\, 1 \,/\, 2}\,f_{\,i},\, f \,|\, a_{\,2},\, \cdots,\, a_{\,n}\,\right>. 
\end{align*}
Thus, \,$\left\{\,S_{C}^{\,-\, 1 \,/\, 2}\,f_{\,i}\,\right\}^{\,\infty}_{\,i \,=\, 1}$\, is a \,$C$-controlled Parseval frame associated to \,$\left(\,a_{\,2},\, \cdots,\, a_{\,n}\,\right)$\, for \,$H$.    
\end{proof}

We now give the concept of dual controlled frame in \,$n$-Hilbert space.

\begin{definition}\label{defi1} 
Let \,$\left\{\,f_{\,i}\,\right\}^{\,\infty}_{i \,=\, 1}$\, be a \,$C$-controlled frame associated to \,$\left(\,a_{\,2} \,,\, \cdots \,,\, a_{\,n}\,\right)$\, for \,$H$.\,Then a \,$C$-controlled frame \,$\left\{\,g_{\,i}\,\right\}^{\,\infty}_{i \,=\, 1}$\, associated to \,$\left(\,a_{\,2} \,,\, \cdots \,,\, a_{\,n}\,\right)$\, satisfying
\[f \,=\, \sum\limits^{\,\infty}_{i \,=\, 1}\, \left<\,f,\, g_{\,i} \,|\, a_{\,2},\, \cdots,\, a_{\,n}\,\right>\,C\,f_{\,i}\; \;\;\forall\; f \,\in\, H_{F}\]
is called a dual controlled frame or alternative dual controlled frame associated to \,$\left(a_{\,2},\, \cdots,\, a_{\,n}\right)$\, of \,$\left\{\,f_{\,i}\,\right\}^{\,\infty}_{i \,=\, 1}$.
\end{definition} 
 
\begin{theorem}\label{th3.1}
Let \,$\left\{\,f_{\,i}\,\right\}^{\,\infty}_{i \,=\, 1}$\, and \,$\left\{\,g_{\,i}\,\right\}^{\,\infty}_{i \,=\, 1}$\, be two \,$C$-controlled Bessel sequences associated to \,$\left(\,a_{\,2} \,,\, \cdots \,,\, a_{\,n}\,\right)$\, in \,$H$.\,Then the following are equivalent:
\begin{itemize}
\item[$(i)$] $f \,=\, \sum\limits^{\,\infty}_{i \,=\, 1}\, \left<\,f,\, g_{\,i} \,|\, a_{\,2},\, \cdots,\, a_{\,n}\,\right>\,C\,f_{\,i}\; \;\forall\; f \,\in\, H_{F}$.
\item[$(ii)$] $f \,=\, \sum\limits^{\,\infty}_{i \,=\, 1}\, \left<\,f,\, g_{\,i} \,|\, a_{\,2},\, \cdots,\, a_{\,n}\,\right>\,C\,g_{\,i}\; \;\forall\; f \,\in\, H_{F}$.
\end{itemize}  
\end{theorem}

\begin{proof}$(\,i\,) \,\Rightarrow\, (\,ii\,)$
Let \,$T_{C}\; \;\text{and}\; \,T_{C^{\,\prime}}$\, be the pre-frame operators of \,$\left\{\,f_{\,i}\,\right\}^{\,\infty}_{i \,=\, 1}$\, and \,$\left\{\,g_{\,i}\,\right\}^{\,\infty}_{i \,=\, 1}$, respectively.\;Composing \,$T_{C}$\, with the adjoint of \,$T_{C^{\,\prime}}$, for all \,$f \,\in\, H_{F}$, we get 
\[T_{C}\,T_{C^{\,\prime}}^{\,\ast} \,:\, H_{F} \,\to\, H_{F},\; T_{C}\,T_{C^{\,\prime}}^{\,\ast}\,f \,=\, \sum\limits^{\,\infty}_{i \,=\, 1}\, \left<\,f,\, g_{\,i} \,|\, a_{\,2},\, \cdots,\, a_{\,n}\,\right>\,C\,f_{\,i}.\]
Now, in terms of pre-frame operators \,$(\,i\,)$\, can be written as \;$T_{C}\,T_{C^{\,\prime}}^{\,\ast} \,=\, I_{F}$\; and this equivalent to \,$T_{C^{\,\prime}}\,T^{\,\ast}_{C} \,=\, I_{F}$\,.\;Therefore, for each \,$f \,\in\, H_{F}$, 
\[f \,=\, \sum\limits^{\,\infty}_{i \,=\, 1}\, \left<\,f,\, f_{\,i} \,|\, a_{\,2},\, \cdots,\, a_{\,n}\,\right>\,C\,g_{\,i}.\]
Similarly, \,$(\,ii\,) \,\Rightarrow\, (\,i\,)$\, follows. 
\end{proof}

\begin{remark}
Suppose that the equivalent conditions of Theorem \ref{th3.1} are satisfied.\;Then using Cauchy-Schwartz inequality, for every \,$f \,\in\, H_{F}$, we have
\begin{align*}
&\left\|\,f,\, a_{\,2},\, \cdots,\, a_{\,n}\,\right\|^{\,2} \,=\, \left<\,f,\, f \,|\, a_{\,2},\, \cdots,\, a_{\,n}\,\right>\\
&\,=\, \left<\,\sum\limits^{\,\infty}_{i \,=\, 1}\, \left<\,f,\, f_{\,i} \,|\, a_{\,2},\, \cdots,\, a_{\,n}\,\right>\,C\,g_{\,i},\, f \,|\, a_{\,2},\, \cdots,\, a_{\,n}\,\right>\\
&=\, \sum\limits^{\,\infty}_{i \,=\, 1}\, \left<\,f,\,f_{\,i} \,|\, a_{\,2},\, \cdots,\, a_{\,n}\,\right>\,\left<\,C\,g_{\,i},\,f \,|\, a_{\,2},\, \cdots,\, a_{\,n}\,\right>\\
& \leq \left(\sum\limits^{\,\infty}_{i \,=\, 1}\left|\,\left<\,f,\, f_{\,i}\,|\,a_{2},\, \cdots,\, a_{n}\,\right>\,\right|^{\,2}\right)^{1 \,/\, 2} \left(\,\sum\limits^{\,\infty}_{i \,=\, 1}\left|\,\left<\,C\,f,\, g_{\,i} \,|\, a_{2},\, \cdots,\, a_{n}\,\right>\,\right|^{\,2}\right)^{1 \,/\, 2}\\
&\leq\, \left(\,\left\|\,C^{\,-\, 1 \,/\, 2}\,\right\|^{\,2}\,\sum\limits_{i \,=\, 1}^{\,\infty}\,\left<\,f,\, f_{\,i} \,|\, a_{\,2},\, \cdots,\, a_{\,n}\,\right>\,\left<\,C\,f_{\,i},\, f \,|\, a_{\,2},\, \cdots,\, a_{\,n}\,\right>\,\right)^{1 \,/\, 2}\,\times\\
& \left(\,\left\|\,C^{\,-\, 1 \,/\, 2}\,\right\|^{\,2}\,\sum\limits_{i \,=\, 1}^{\,\infty}\,\left<\,C\,f,\, g_{\,i} \,|\, a_{\,2},\, \cdots,\, a_{\,n}\,\right>\,\left<\,C\,g_{\,i},\, C\,f \,|\, a_{\,2},\, \cdots,\, a_{\,n}\,\right>\,\right)^{1 \,/\, 2}\\
&\leq\,\left(\,\sum\limits_{i \,=\, 1}^{\,\infty}\,\left<\,f,\, f_{\,i} \,|\, a_{\,2},\, \cdots,\, a_{\,n}\,\right>\,\left<\,C\,f_{\,i},\, f \,|\, a_{\,2},\, \cdots,\, a_{\,n}\,\right>\,\right)^{1 \,/\, 2}\,\times\\
&\hspace{2cm} B^{1 \,/\, 2}\,\left\|\,C^{\,-\, 1 \,/\, 2}\,\right\|^{\,2}\, \left\|\,C\,f,\, a_{\,2},\, \cdots,\, a_{\,n}\,\right\|_{1}\\
&\Rightarrow\, \dfrac{1}{B\,\left\|\,C^{\,-\, 1 \,/\, 2}\,\right\|^{\,4}\,\|\,C\,\|^{\,2}}\; \left\|\,f,\, a_{\,2},\, \cdots,\, a_{\,n}\,\right\|_{1}^{\,2}\\
&\hspace{1cm} \,\leq\, \sum\limits_{i \,=\, 1}^{\,\infty}\,\left<\,f,\, f_{\,i} \,|\, a_{\,2},\, \cdots,\, a_{\,n}\,\right>\,\left<\,C\,f_{\,i},\, f \,|\, a_{\,2},\, \cdots,\, a_{\,n}\,\right>.
\end{align*}
This shows that \,$\left\{\,f_{\,i}\,\right\}^{\,\infty}_{i \,=\, 1}$\, is a \,$C$-controlled frame associated to \,$\left(\,a_{\,2},\, \cdots,\, a_{\,n}\,\right)$\, for \,$H$.\;Similarly, it can be shown that \,$\left\{\,g_{\,i}\,\right\}^{\,\infty}_{i \,=\, 1}$\, is also a \,$C$-controlled frame associated to \,$\left(\,a_{\,2},\, \cdots,\, a_{\,n}\,\right)$\, for \,$H$.
\end{remark}

\section{Controlled frame in tensor product of $n$-Hilbert spaces}

\smallskip\hspace{.6 cm}In this section, we present the concept of controlled frame in tensor product of \,$n$-Hilbert spaces and give a characterization.\,Dual controlled frame in tensor product of \,$n$-Hilbert spaces is also described.\,Finally, we consider the direct sum of controlled frames in \,$n$-Hilbert spaces.\\

Let \,$H$\, and \,$K$\, be two \,$n$-Hilbert spaces associated with the \,$n$-inner products \,$\left<\,\cdot,\, \cdot \,|\, \cdot,\, \cdots,\, \cdot\,\right>_{1}$\, and \,$\left<\,\cdot,\, \cdot \,|\, \cdot,\, \cdots,\, \cdot\,\right>_{2}$, respectively.\;The tensor product of \,$H$\, and \,$K$\, is denoted by \,$H \,\otimes\, K$\, and it is defined to be an \,$n$-inner product space associated with the \,$n$-inner product given by 
\[\left<\,f_{\,1} \,\otimes\, g_{\,1},\, f_{\,2} \,\otimes\, g_{\,2} \,|\, f_{\,3} \,\otimes\, g_{\,3},\, \,\cdots,\, f_{\,n} \,\otimes\, g_{\,n}\,\right>\]
\begin{equation}\label{eqn1}
 \,=\, \left<\,f_{\,1},\, f_{\,2} \,|\, f_{\,3},\, \,\cdots,\, f_{\,n}\,\right>_{1}\,\left<\,g_{\,1},\, g_{\,2} \,|\, g_{\,3},\, \,\cdots,\, g_{\,n}\,\right>_{2},
\end{equation}
for all \,$f_{\,1},\, f_{\,2},\, f_{\,3},\, \,\cdots,\, f_{\,n} \,\in\, H$\, and \,$g_{\,1},\, g_{\,2},\, g_{\,3},\, \,\cdots,\, g_{\,n} \,\in\, K$.\\
The \,$n$-norm on \,$H \,\otimes\, K$\, is defined by 
\[\left\|\,f_{\,1} \,\otimes\, g_{\,1},\, f_{\,2} \,\otimes\, g_{\,2},\, \,\cdots,\,\, f_{\,n} \,\otimes\, g_{\,n}\,\right\|\]
\begin{equation}\label{eqn1.1}
\hspace{.6cm} =\,\left\|\,f_{\,1},\, f_{\,2},\, \cdots,\, f_{\,n}\,\right\|_{1}\;\left\|\,g_{\,1},\, g_{\,2},\, \cdots,\, g_{\,n}\,\right\|_{2},
\end{equation}
for all \,$f_{\,1},\, f_{\,2},\, \,\cdots,\, f_{\,n} \,\in\, H\, \;\text{and}\; \,g_{\,1},\, g_{\,2},\, \,\cdots,\, g_{\,n} \,\in\, K$, where the \,$n$-norms \,$\left\|\,\cdot,\, \cdots,\, \cdot \,\right\|_{1}$\, and \,$\left\|\,\cdot,\, \cdots,\, \cdot \,\right\|_{2}$\, are generated by \,$\left<\,\cdot,\, \cdot \,|\, \cdot,\, \cdots,\, \cdot\,\right>_{1}$\, and \,$\left<\,\cdot,\, \cdot \,|\, \cdot,\, \cdots,\, \cdot\,\right>_{2}$, respectively.\;The space \,$H \,\otimes\, K$\, is complete with respect to the above \,$n$-inner product.\;Therefore the space \,$H \,\otimes\, K$\, is an \,$n$-Hilbert space.\\

Consider \,$G \,=\, \left\{\,b_{\,2},\, b_{\,3},\, \cdots,\, b_{\,n}\,\right\}$, where \,$b_{\,2},\, b_{\,3},\, \cdots,\, b_{\,n}$\, are fixed elements in \,$K$\, and \,$L_{G}$\, denote the linear subspace of \,$K$\, spanned by \,$G$.\,Now, we can define the Hilbert space \,$K_{G}$\, with respect to the inner product is given by
\[\left<\,p \,+\, L_{G}\,,\, q \,+\, L_{G}\,\right>_{G} \,=\, \left<\,p \,,\, q\,\right>_{G} \,=\, \left<\,p \,,\, q \,|\, b_{\,2} \,,\,  \cdots \,,\, b_{\,n} \,\right>_{2}; \;\forall \;\; p,\, q \,\in\, K.\]

\begin{remark}
According to the definition \ref{def0.001}, \,$H_{F} \,\otimes\, K_{G}$\, is the Hilbert space with respect to the inner product:
\[\left<\,p \,\otimes\, q \,,\, p^{\,\prime} \,\otimes\, q^{\,\prime}\,\right> \,=\, \left<\,p \,,\, p^{\,\prime}\,\right>_{F}\;\left<\,q \,,\, q^{\,\prime}\,\right>_{G},\]
for all \,$p,\, p^{\,\prime} \,\in\, H_{F}\; \;\text{and}\; \;q,\, q^{\,\prime} \,\in\, K_{G}$.    
\end{remark}

\begin{definition}
Let \,$C_{1} \,\otimes\, C_{2} \,\in\, \mathcal{G}\,\mathcal{B}\left(\,H_{F} \,\otimes\, K_{G}\,\right)$.\,Then the sequence \,$\left\{\,f_{\,i} \,\otimes\, g_{\,j}\,\right\}^{\,\infty}_{i,\,j \,=\, 1}$\,  in \,$H \,\otimes\, K$\, is said to be a \,$\left(\,C_{1} \,\otimes\, C_{2}\,\right)$-controlled frame associated to \,$(\,a_{\,2} \,\otimes\, b_{\,2},\, \,\cdots,\, a_{\,n} \,\otimes\, b_{\,n}\,)$\, for \,$H \,\otimes\, K$\, if there exist constants \,$0 \,<\, A \,\leq\, B \,<\, \infty$\, such that
\begin{align}
&A \left\|\,f \,\otimes\, g,\, a_{\,2} \,\otimes\, b_{\,2},\,\cdots,\, a_{\,n} \,\otimes\, b_{\,n}\,\right\|^{\,2}\label{eqp1.1}\\
& \leq \sum\limits_{i,\, j \,=\, 1}^{\,\infty}\,\left<\,f \,\otimes\, g,\, f_{\,i} \,\otimes\, g_{\,j} \,|\, a_{\,2} \,\otimes\, b_{\,2},\, \cdots,\, a_{\,n} \,\otimes\, b_{\,n}\,\right>\,\times\nonumber\\
&\hspace{1.5cm}\left<\,\left(\,C_{1} \,\otimes\, C_{2}\,\right)\,\left(\,f_{\,i} \,\otimes\, g_{\,j}\,\right),\, f \,\otimes\, g \,|\, a_{\,2} \,\otimes\, b_{\,2},\, \cdots,\, a_{\,n} \,\otimes\, b_{\,n}\,\right>\nonumber\\
&\leq\, B\,\left\|\,f \,\otimes\, g,\, a_{\,2} \,\otimes\, b_{\,2},\, \cdots,\, a_{\,n} \,\otimes\, b_{\,n}\,\right\|^{\,2}; \;\forall\, \;f \,\otimes\, g \,\in\, H_{F} \,\otimes\, K_{G},\nonumber
\end{align}
where \,$\{\,f_{\,i}\,\}_{i \,=\,1}^{\infty}$\, and \,$\{\,g_{\,j}\,\}_{j \,=\,1}^{\infty}$\, be the sequences of vectors in \,$H$\, and \,$K$, respectively and \,$a_{\,2} \,\otimes\, b_{\,2},\, a_{\,3} \,\otimes\, b_{\,3},\,\cdots,\, a_{\,n} \,\otimes\, b_{\,n}$\, be fixed elements in \,$H \,\otimes\, K$.\;The constants \,$A,\,B$\, are called the frame bounds.\,If the sequence \,$\left\{\,f_{\,i} \,\otimes\, g_{\,j}\,\right\}^{\,\infty}_{i,\,j \,=\, 1}$\, satisfies only the right inequality of (\ref{eqp1.1}) then it is called a \,$\left(\,C_{1} \,\otimes\, C_{2}\,\right)$-controlled Bessel sequence associated to \,$\left(\,a_{\,2} \,\otimes\, b_{\,2},\, \,\cdots,\, a_{\,n} \,\otimes\, b_{\,n}\,\right)$\, in \,$H \,\otimes\, K$.    
\end{definition}

\begin{theorem}\label{th2.1}
Let \,$\{\,f_{\,i}\,\}_{i \,=\,1}^{\infty}$\, and \,$\{\,g_{\,j}\,\}_{j \,=\,1}^{\infty}$\, be the sequences of vectors in \,$n$-Hilbert spaces \,$H$\, and \,$K$.\,The sequence \,$\left\{\,f_{\,i} \,\otimes\, g_{\,j}\,\right\}^{\,\infty}_{i,\,j \,=\, 1} \,\subseteq\, H \,\otimes\, K$\, is a \,$\left(\,C_{1} \,\otimes\, C_{2}\,\right)$-controlled frame associated to \,$\left(\,a_{\,2} \,\otimes\, b_{\,2},\, \,\cdots,\, a_{\,n} \,\otimes\, b_{\,n}\,\right)$\, for \,$H \,\otimes\, K$\, if and only if \,$\{\,f_{\,i}\,\}_{i \,=\,1}^{\infty}$\, is a \,$C_{1}$-controlled frame associated to \,$\left(\,a_{\,2},\, \cdots,\, a_{\,n}\,\right)$\, for \,$H$\, and \,$\{\,g_{\,j}\,\}_{j \,=\,1}^{\infty}$\, is a \,$C_{2}$-controlled frame associated to \,$\left(\,b_{\,2},\, \cdots,\, b_{\,n}\,\right)$\, for \,$K$.   
\end{theorem}

\begin{proof}
Suppose that the sequence \,$\left\{\,f_{\,i} \,\otimes\, g_{\,j}\,\right\}^{\,\infty}_{i,\,j \,=\, 1}$\, is a \,$\left(\,C_{1} \,\otimes\, C_{2}\,\right)$-controlled frame associated to \,$\left(\,a_{\,2} \,\otimes\, b_{\,2},\, \,\cdots,\, a_{\,n} \,\otimes\, b_{\,n}\,\right)$\, for \,$H \,\otimes\, K$.\;Then, for each \,$f \,\otimes\, g \,\in\, H_{F} \,\otimes\, K_{G} \,-\, \{\,\theta \,\otimes\, \theta\,\}$, there exist constants \,$A,\,B \,>\, 0$\, such that
\begin{align*}
&A \left\|\,f \,\otimes\, g,\, a_{\,2} \,\otimes\, b_{\,2},\,\cdots,\, a_{\,n} \,\otimes\, b_{\,n}\,\right\|^{\,2}\\
& \leq \sum\limits_{i,\, j \,=\, 1}^{\,\infty}\,\left<\,f \,\otimes\, g,\, f_{\,i} \,\otimes\, g_{\,j} \,|\, a_{\,2} \,\otimes\, b_{\,2},\, \cdots,\, a_{\,n} \,\otimes\, b_{\,n}\,\right>\,\times\\
&\hspace{1.5cm}\left<\,\left(\,C_{1} \,\otimes\, C_{2}\,\right)\,\left(\,f_{\,i} \,\otimes\, g_{\,j}\,\right),\, f \,\otimes\, g \,|\, a_{\,2} \,\otimes\, b_{\,2},\, \cdots,\, a_{\,n} \,\otimes\, b_{\,n}\,\right>\\
&\leq\, B\,\left\|\,f \,\otimes\, g,\, a_{\,2} \,\otimes\, b_{\,2},\, \cdots,\, a_{\,n} \,\otimes\, b_{\,n}\,\right\|^{\,2}
\end{align*}
Using (\ref{eqn1}) and (\ref{eqn1.1}), we get
\begin{align*}
&A\left\|\,f,\, a_{\,2},\, \cdots,\, a_{\,n}\,\right\|_{1}^{\,2}\,\left\|\,g,\, b_{\,2},\, \cdots,\, b_{\,n}\,\right\|_{2}^{\,2}\\
&\leq\,\sum\limits_{i,\, j \,=\, 1}^{\,\infty}\,\left<\,f,\, f_{\,i}\,|\,a_{\,2},\, \cdots,\, a_{\,n}\,\right>_{1}\,\left<\,g,\, g_{\,j} \,|\, b_{\,2},\, \cdots,\, b_{\,n}\,\right>_{2}\,\times\\
&\hspace{1.5cm}\left<\,C_{1}\,f_{\,i},\, f \,|\,a_{\,2},\, \cdots,\, a_{\,n}\,\right>_{1}\,\left<\,C_{2}\,g_{\,j},\, g \,|\, b_{\,2},\, \cdots,\, b_{\,n}\,\right>_{2}\\
&\leq\,B\,\left\|\,f,\, a_{\,2},\, \cdots,\, a_{\,n}\,\right\|_{1}^{\,2}\,\left\|\,g,\, b_{\,2},\, \cdots,\, b_{\,n}\,\right\|_{2}^{\,2}\\
&\Rightarrow\, A\left\|\,f,\, a_{\,2},\, \cdots,\, a_{\,n}\,\right\|_{1}^{\,2}\,\left\|\,g,\, b_{\,2},\, \cdots,\, b_{\,n}\,\right\|_{2}^{\,2}\\
&\sum\limits_{\,i \,=\, 1}^{\,\infty}\,\left<\,f ,\, f_{\,i} \,|\, a_{\,2},\, \cdots,\, a_{\,n}\,\right>_{1}\,\left<\,C_{1}\,f_{\,i},\, f \,|\,a_{\,2},\, \cdots,\, a_{\,n}\,\right>_{1}\,\times\\
&\hspace{1.5cm}\sum\limits_{\,j \,=\, 1}^{\,\infty}\,\left<\,g,\, g_{\,j} \,|\, b_{\,2},\, \cdots,\, b_{\,n}\,\right>_{2}\,\left<\,C_{2}\,g_{\,j},\, g \,|\, b_{\,2},\, \cdots,\, b_{\,n}\,\right>_{2}\\
&\leq\,B\,\left\|\,f,\, a_{\,2},\, \cdots,\, a_{\,n}\,\right\|_{1}^{\,2}\,\left\|\,g,\, b_{\,2},\, \cdots,\, b_{\,n}\,\right\|_{2}^{\,2}
\end{align*}
Since \,$f \,\otimes\, g \,\in\, H_{F} \,\otimes\, K_{G}$\, is  non-zero element i.\,e., \,$f \,\in\, H_{F}$\, and \,$g \,\in\, K_{G}$\, are non-zero elements.\,Here, we may assume that every \,$f_{\,i},\, C_{1}\,f_{\,i}$\, and \,$a_{\,2},\, \cdots$, \,$ a_{\,n}$\, are linearly independent and every \,$g_{\,j},\, C_{2}\,g_{\,j}$\, and \,$b_{\,2},\, \cdots$, \,$ b_{\,n}$\, are linearly independent.\,Hence 
\[\sum\limits_{\,j \,=\, 1}^{\,\infty}\,\left<\,g,\, g_{\,j} \,|\, b_{\,2},\, \cdots,\, b_{\,n}\,\right>_{2}\,\left<\,C_{2}\,g_{\,j},\, g \,|\, b_{\,2},\, \cdots,\, b_{\,n}\,\right>_{2},\]
\[\sum\limits_{\,i \,=\, 1}^{\,\infty}\,\left<\,f ,\, f_{\,i} \,|\, a_{\,2},\, \cdots,\, a_{\,n}\,\right>_{1}\,\left<\,C_{1}\,f_{\,i},\, f \,|\,a_{\,2},\, \cdots,\, a_{\,n}\,\right>_{1}\] are non-zero.\,Therefore, by the above inequality, we get
\begin{align*}
&\dfrac{A\,\left\|\,g,\, b_{\,2},\, \cdots,\, b_{\,n}\,\right\|_{2}^{\,2}\,\left\|\,f,\, a_{\,2},\, \cdots,\, a_{\,n}\,\right\|_{1}^{\,2}}{\sum\limits_{\,j \,=\, 1}^{\,\infty}\,\left<\,g,\, g_{\,j} \,|\, b_{\,2},\, \cdots,\, b_{\,n}\,\right>_{2}\,\left<\,C_{2}\,g_{\,j},\, g \,|\, b_{\,2},\, \cdots,\, b_{\,n}\,\right>_{2}}\\
&\leq\,\sum\limits_{\,i \,=\, 1}^{\,\infty}\,\left<\,f ,\, f_{\,i} \,|\, a_{\,2},\, \cdots,\, a_{\,n}\,\right>_{1}\,\left<\,C_{1}\,f_{\,i},\, f \,|\,a_{\,2},\, \cdots,\, a_{\,n}\,\right>_{1}\\
&\leq\,\dfrac{B\,\left\|\,g,\, b_{\,2},\, \cdots,\, b_{\,n}\,\right\|_{2}^{\,2}\,\left\|\,f,\, a_{\,2},\, \cdots,\, a_{\,n}\,\right\|_{1}^{\,2}}{\sum\limits_{\,j \,=\, 1}^{\,\infty}\,\left<\,g,\, g_{\,j} \,|\, b_{\,2},\, \cdots,\, b_{\,n}\,\right>_{2}\,\left<\,C_{2}\,g_{\,j},\, g \,|\, b_{\,2},\, \cdots,\, b_{\,n}\,\right>_{2}} 
\end{align*}   
This implies that
\begin{align*}
A_{1} \, \left\|\,f,\, a_{\,2},\, \cdots,\, a_{\,n} \,\right\|_{1}^{\,2}& \,\leq\, \sum\limits_{i \,=\, 1}^{\,\infty}\,\left<\,f,\, f_{\,i} \,|\, a_{\,2},\, \cdots,\, a_{\,n}\,\right>_{1}\,\left<\,C_{1}\,f_{\,i},\, f \,|\, a_{\,2},\, \cdots,\, a_{\,n}\,\right>_{1}\nonumber\\
&\,\leq\, B_{1}\, \left\|\,f,\, a_{\,2},\, \cdots,\, a_{\,n}\,\right\|_{1}^{\,2}\; \;\forall\; f \,\in\, H_{F},
\end{align*}
where \,$A_{1} \,=\, \inf\limits_{g \,\in\, K_{G}}\,\left\{\,\dfrac{A\,\left\|\,g,\, b_{\,2},\, \cdots,\, b_{\,n}\,\right\|_{2}^{\,2}}{\sum\limits_{\,j \,=\, 1}^{\,\infty}\,\left<\,g,\, g_{\,j} \,|\, b_{\,2},\, \cdots,\, b_{\,n}\,\right>_{2}\,\left<\,C_{2}\,g_{\,j},\, g \,|\, b_{\,2},\, \cdots,\, b_{\,n}\,\right>_{2}}\,\right\}$\, and \\$B_{1} \,=\, \sup\limits_{g \,\in\, K_{G}}\,\left\{\,\dfrac{B\,\left\|\,g,\, b_{\,2},\, \cdots,\, b_{\,n}\,\right\|_{2}^{\,2}}{\sum\limits_{\,j \,=\, 1}^{\,\infty}\,\left<\,g,\, g_{\,j} \,|\, b_{\,2},\, \cdots,\, b_{\,n}\,\right>_{2}\,\left<\,C_{2}\,g_{\,j},\, g \,|\, b_{\,2},\, \cdots,\, b_{\,n}\,\right>_{2}}\,\right\}$.\\This shows that \,$\{\,f_{\,i}\,\}_{i \,=\,1}^{\infty}$\, is a \,$C_{1}$-controlled frame associated to \,$\left(\,a_{\,2},\, \cdots,\, a_{\,n}\,\right)$\, for \,$H$.\;Similarly, it can be shown that \,$\{\,g_{\,j}\,\}_{j \,=\,1}^{\infty}$\, is a \,$C_{2}$-controlled frame associated to \,$\left(\,b_{\,2},\, \cdots,\, b_{\,n}\,\right)$\, for \,$K$.\\

Conversely, suppose that \,$\{\,f_{\,i}\,\}_{i \,=\,1}^{\infty}$\, is a \,$C_{1}$-controlled frame associated to \,$(\,a_{2},\, \cdots,\, a_{n}\,)$\, for \,$H$\, with bounds \,$A,\,B$\, and \,$\{\,g_{\,j}\,\}_{j \,=\,1}^{\infty}$\, is a \,$C_{2}$-controlled frame associated to \,$\left(\,b_{2},\, \cdots,\, b_{n}\,\right)$\, for \,$K$\, with bounds \,$C,\,D$.\,Then
\begin{align*}
A \, \left\|\,f,\, a_{\,2},\, \cdots,\, a_{\,n} \,\right\|_{1}^{\,2}& \,\leq\, \sum\limits_{i \,=\, 1}^{\,\infty}\,\left<\,f,\, f_{\,i} \,|\, a_{\,2},\, \cdots,\, a_{\,n}\,\right>_{1}\,\left<\,C_{1}\,f_{\,i},\, f \,|\, a_{\,2},\, \cdots,\, a_{\,n}\,\right>_{1}\nonumber\\
&\,\leq\, B\, \left\|\,f,\, a_{\,2},\, \cdots,\, a_{\,n}\,\right\|_{1}^{\,2}\; \;\forall\; f \,\in\, H_{F},
\end{align*}
\begin{align*}
C\left\|\,g,\, b_{\,2},\, \cdots,\, b_{\,n}\,\right\|_{2}^{\,2} \,&\leq\,\sum\limits_{\,j \,=\, 1}^{\,\infty}\,\left<\,g,\, g_{\,j} \,|\, b_{\,2},\, \cdots,\, b_{\,n}\,\right>_{2}\,\left<\,C_{2}\,g_{\,j},\, g \,|\, b_{\,2},\, \cdots,\, b_{\,n}\,\right>_{2}\\
&\leq\, D\,\left\|\,g,\, b_{\,2},\, \cdots,\, b_{\,n}\,\right\|_{2}^{\,2}\; \;\forall\; g \,\in\, K_{G}. 
\end{align*}
Multiplying the above two inequalities and using (\ref{eqn1}) and (\ref{eqn1.1}), we get
\begin{align*}
&A\,C\,\left\|\,f \,\otimes\, g,\, a_{\,2} \,\otimes\, b_{\,2},\,\cdots,\, a_{\,n} \,\otimes\, b_{\,n}\,\right\|^{\,2}\\
& \leq \sum\limits_{i,\, j \,=\, 1}^{\,\infty}\,\left<\,f \,\otimes\, g,\, f_{\,i} \,\otimes\, g_{\,j} \,|\, a_{\,2} \,\otimes\, b_{\,2},\, \cdots,\, a_{\,n} \,\otimes\, b_{\,n}\,\right>\,\times\\
&\hspace{1.5cm}\left<\,\left(\,C_{1} \,\otimes\, C_{2}\,\right)\,\left(\,f_{\,i} \,\otimes\, g_{\,j}\,\right),\, f \,\otimes\, g \,|\, a_{\,2} \,\otimes\, b_{\,2},\, \cdots,\, a_{\,n} \,\otimes\, b_{\,n}\,\right>\\
&\leq\, B\,D\,\left\|\,f \,\otimes\, g,\, a_{\,2} \,\otimes\, b_{\,2},\, \cdots,\, a_{\,n} \,\otimes\, b_{\,n}\,\right\|^{\,2}; \;\forall\, \;f \,\otimes\, g \,\in\, H_{F} \,\otimes\, K_{G}.
\end{align*}
Hence, \,$\left\{\,f_{\,i} \,\otimes\, g_{\,j}\,\right\}^{\,\infty}_{i,\,j \,=\, 1}$\, is a \,$\left(\,C_{1} \,\otimes\, C_{2}\,\right)$-controlled frame associated to \,$(\,a_{\,2} \otimes b_{\,2},\, \,\cdots,\, a_{\,n} \otimes b_{\,n}\,)$\, for \,$H \,\otimes\, K$.\,This completes the proof. 
\end{proof}

\begin{remark}
Let \,$\left\{\,f_{\,i} \,\otimes\, g_{\,j}\,\right\}^{\,\infty}_{i,\,j \,=\, 1}$\, be a \,$\left(\,C_{1} \,\otimes\, C_{2}\,\right)$-controlled frame associated to \,$(\,a_{\,2} \,\otimes\, b_{\,2}$, \,$\cdots,\, a_{\,n} \,\otimes\, b_{\,n}\,)$\, for \,$H \,\otimes\, K$.\;According to the definition \ref{def0.0001}, the frame operator \,$S_{C_{1} \,\otimes\, C_{2}} \,:\, H_{F} \,\otimes\, K_{G} \,\to\, H_{F} \,\otimes\, K_{G}$\, is described by
\begin{align*}
&S_{C_{1} \,\otimes\, C_{2}}\,(\,f \,\otimes\, g\,)\\
& = \sum\limits_{i,\, j \,=\, 1}^{\,\infty}\left<\,f \otimes g,\, f_{\,i} \otimes g_{\,j} \,|\, a_{\,2} \otimes b_{\,2}, \,\cdots,\, a_{\,n} \otimes b_{\,n}\,\right>\,\left(\,C_{1} \,\otimes\, C_{2}\,\right)\left(\,f_{\,i} \otimes g_{\,j}\,\right)
\end{align*} 
for all \,$f \,\otimes\, g \,\in\, H_{F} \,\otimes\, K_{G}$. 
\end{remark}

\begin{proposition}
If \,$S_{C_{1}},\, \,S_{C_{2}}$\, and \,$S_{C_{1} \,\otimes\, C_{2}}$\, are the corresponding frame operator for \,$\{\,f_{\,i}\,\}_{i \,=\,1}^{\infty},\, \,\{\,g_{\,j}\,\}_{j \,=\,1}^{\infty}$\, and \,$\left\{\,f_{\,i} \,\otimes\, g_{\,j}\,\right\}^{\,\infty}_{i,\,j \,=\, 1}$, respectively, then \,$S_{C_{1} \,\otimes\, C_{2}} \,=\, S_{C_{1}} \,\otimes\, S_{C_{2}}$\, and \,$S_{C_{1} \,\otimes\, C_{2}}^{\,-\, 1} \,=\, S^{\,-\, 1}_{C_{1}} \,\otimes\, S^{\,-\, 1}_{C_{2}}$.  
\end{proposition}

\begin{proof}
Since \,$S_{C_{1} \,\otimes\, C_{2}}$\, is the frame operator for \,$\left\{\,f_{\,i} \,\otimes\, g_{\,j}\,\right\}^{\,\infty}_{i,\,j \,=\, 1}$, we have
\begin{align*}
&S_{C_{1} \,\otimes\, C_{2}}\,(\,f \,\otimes\, g\,)\\
& = \sum\limits_{i,\, j \,=\, 1}^{\,\infty}\left<\,f \otimes g,\, f_{\,i} \otimes g_{\,j} \,|\, a_{\,2} \otimes b_{\,2}, \,\cdots,\, a_{\,n} \otimes b_{\,n}\,\right>\,\left(\,C_{1} \,\otimes\, C_{2}\,\right)\left(\,f_{\,i} \otimes g_{\,j}\,\right)\\
& \,=\, \sum\limits_{i,\, j \,=\, 1}^{\,\infty}\,\left<\,f,\, f_{\,i} \,|\, a_{\,2},\, \cdots,\, a_{\,n}\,\right>_{1}\,\left<\,g,\, g_{\,j} \,|\, b_{\,2},\, \cdots,\, b_{\,n}\,\right>_{2}\,\left(\,C_{1}\,f_{\,i} \,\otimes\, C_{2}\,g_{\,j}\,\right)\\
&= \left(\sum\limits_{\,i \,=\, 1}^{\,\infty}\left<\,f,\, f_{\,i} \,|\, a_{\,2},\, \cdots,\, a_{\,n}\,\right>_{1}\,C_{1}\,f_{\,i}\right) \,\otimes\, \left(\sum\limits_{\,j \,=\, 1}^{\,\infty}\,\left<\,g,\, g_{\,j} \,|\, b_{\,2},\, \cdots,\, b_{\,n}\,\right>_{2}\,C_{2}\,g_{\,j}\right)\\
&=\, S_{C_{1}}\,f \,\otimes\, S_{C_{2}}\,g \,=\, \left(\,S_{C_{1}} \,\otimes\, S_{C_{2}}\,\right)\,(\,f \,\otimes\, g\,)\; \;\forall\; f \,\otimes\, g \,\in\, H_{F} \,\otimes\, K_{G}.
\end{align*}
Thus, \,$S_{C_{1} \,\otimes\, C_{2}} \,=\, S_{C_{1}} \,\otimes\, S_{C_{2}}$.\,Since \,$S_{C_{1}}$\, and \,$S_{C_{2}}$\, are invertible, by Theorem \ref{th1.1} \,$(iv)$, \,$S_{C_{1} \,\otimes\, C_{2}}^{\,-\, 1} \,=\, \left(\,S_{C_{1}} \,\otimes\, S_{C_{2}}\,\right)^{\,-\, 1} \,=\, S^{\,-\, 1}_{C_{1}} \,\otimes\, S^{\,-\, 1}_{C_{2}}$.\,This completes the proof.    
\end{proof}

\begin{proposition}
Let \,$\{\,f_{\,i}\,\}_{i \,=\,1}^{\infty}$\, be a \,$C_{1}$-controlled frame associated to \,$(\,a_{\,2},\, \cdots,\\ a_{\,n}\,)$\, for \,$H$\, with bounds \,$A,\,B$ and \,$\{\,g_{\,j}\,\}_{j \,=\,1}^{\infty}$\, be a \,$C_{2}$-controlled frame associated to \,$\left(\,b_{\,2},\, \cdots,\, b_{\,n}\,\right)$\, for \,$K$\, with bounds \,$C,\,D$\, having their corresponding frame operator \,$S_{C_{1}}$\, and \,$S_{C_{2}}$, respectively.\,Then \,$A\,C\,I_{F \,\otimes\, G} \,\leq\, S_{C_{1} \,\otimes\, C_{2}} \,\leq\, B\,D\,I_{F \,\otimes\, G}$, where \,$I_{F \,\otimes\, G}$\, is the identity operator on \,$H_{F} \,\otimes\, K_{G}$.
\end{proposition}

\begin{proof}
Since \,$S_{C_{1}}$\, and \,$S_{C_{2}}$\, are frame operators, we have 
\[A\,I_{F} \,\leq\, S_{C_{1}} \,\leq\, B\,I_{F},\; \;C\,I_{G} \,\leq\, S_{C_{2}} \,\leq\, D\,I_{G},\]
where \,$I_{F}$\,and \,$I_{G}$\, are the identity operators on \,$H_{F}$\, and \,$K_{G}$, respectively.\,Taking tensor product on the above two inequalities, we get
\begin{align*}
&A\,C \left(\,I_{F} \,\otimes\, I_{G}\,\right) \,\leq\,  \left(\,S_{C_{1}} \,\otimes\, S_{C_{2}}\,\right) \,\leq\, B\,D\,\left(\,I_{F} \,\otimes\, I_{G}\,\right)\\
&\Rightarrow\,A\,C\,I_{F \,\otimes\, G} \,\leq\, S_{C_{1} \,\otimes\, C_{2}} \,\leq\, B\,D\,I_{F \,\otimes\, G}. 
\end{align*}
This completes the proof.
\end{proof}

Next, we establish the frame decomposition formula in \,$H \,\otimes\, K$.

\begin{proposition}\label{eq2.4}
Let \,$\left\{\,f_{\,i}\,\right\}^{\infty}_{i \,=1\, }$\, and \,$\left\{\,g_{\,j}\,\right\}^{\infty}_{j \,=1\, }$\, be \,$C_{1}$-controlled and \,$C_{2}$-controlled frames associated to \,$\left(\,a_{\,2},\, \cdots,\, a_{\,n}\,\right)$\, and \,$\left(\,b_{\,2},\, \cdots,\, b_{\,n}\,\right)$\, for \,$H$\, and \,$K$\, with the corresponding frame operators \,$S_{C_{1}}$\, and \,$S_{C_{2}}$, respectively.\,Then for each \,$f \,\otimes\, g \,\in\, H_{F} \,\otimes\, K_{G}$, we have
\begin{align*}
& f \,\otimes\, g \\
 &\,=\, \sum\limits^{\infty}_{i,\,j \,=\, 1}\,\left <\,f \,\otimes\, g,\, S^{\,-1}_{C_{1} \,\otimes\, C_{2}}\,\left(\,f_{\,i} \,\otimes\, g_{\,j}\,\right) \,|\, a_{\,2} \otimes b_{\,2}, \,\cdots,\, a_{\,n} \otimes b_{\,n}\,\right >\,\left(\,C_{1} \,\otimes\, C_{2}\,\right)\,\left(\,f_{\,i} \,\otimes\, g_{\,j}\,\right), \\
& f \,\otimes\, g \\
 &\,=\, \sum\limits^{\infty}_{i,\,j \,=\, 1}\,\left <\,f \,\otimes\, g,\, f_{\,i} \,\otimes\, g_{\,j} \,|\, a_{\,2} \otimes b_{\,2}, \,\cdots,\, a_{\,n} \otimes b_{\,n}\,\right >\,S^{\,-1}_{C_{1} \,\otimes\, C_{2}}\,\left(\,C_{1} \,\otimes\, C_{2}\,\right)\,\left(\,f_{\,i} \,\otimes\, g_{\,j}\,\right).\\  
\end{align*}
provided both the series converges unconditionally for all \,$f \,\otimes\, g \,\in\, H_{F} \,\otimes\, K_{G}$.
\end{proposition}

\begin{proof}
Since \,$S_{C_{1}}$\, is the corresponding frame operator for \,$\left\{\,f_{\,i}\,\right\}^{\infty}_{i \,=1\, }$.\,Then
\begin{align*}
f \,=\, S_{C_{1}}\, S^{\,-1}_{C_{1}}\,f &\,=\, \sum\limits^{\infty}_{i \,=\, 1}\,\left <\,S^{\,-1}_{C_{1}}\,f,\, f_{\,i} \,|\, a_{\,2},\, \cdots,\, a_{\,n}\,\,\right >\,C_{1}\,f_{\,i}\\
& \;=\; \sum\limits^{\infty}_{i \,=\, 1}\,\left <\,f,\, S^{\,-1}_{C_{1}}\,f_{\,i} \,|\, a_{\,2},\, \cdots,\, a_{\,n}\,\,\right >\,C_{1}\,f_{\,i}\; \;\forall\; f \,\in\, H_{F}.\\
\end{align*}
Similarly, it can be shown that
\[g \,=\, \sum\limits^{\infty}_{j \,=\, 1}\,\left <\,g,\, S^{\,-1}_{C_{2}}\,g_{\,j} \,|\, b_{\,2},\, \cdots,\, b_{\,n}\,\,\right >\,C_{2}\,g_{\,j}\; \;\forall\; g \,\in\, K_{G}.\]
Thus, for each \,$f \,\otimes\, g \,\in\, H_{F} \,\otimes\, K_{G}$, we have
\begin{align*}
& f \,\otimes\, g \\
&\,= \left(\sum\limits^{\infty}_{i \,=\, 1}\left <\,f,\, S^{\,-1}_{C_{1}}\,f_{\,i} \,|\, a_{\,2},\, \cdots,\, a_{\,n}\,\,\right>\,C_{1}\,f_{\,i}\right) \otimes \left(\sum\limits^{\infty}_{j \,=\, 1}\left <\,g,\, S^{\,-1}_{C_{2}}\,g_{\,j} \,|\, b_{\,2},\, \cdots,\, b_{\,n}\,\,\right >\,C_{2}\,g_{\,j}\right)\\
&=\,\sum\limits^{\infty}_{i,\,j \,=\, 1}\,\left <\,f,\, S^{\,-1}_{C_{1}}\,f_{\,i} \,|\, a_{\,2},\, \cdots,\, a_{\,n}\,\,\right>_{1}\,\left <\,g,\, S^{\,-1}_{C_{2}}\,g_{\,j} \,|\, b_{\,2},\, \cdots,\, b_{\,n}\,\,\right>_{2}\,\left(\,C_{1}\,f_{\,i} \,\otimes\, C_{2}\,g_{\,j}\,\right)\\
&=\sum\limits^{\infty}_{i,\,j \,=\, 1}\,\left <\,f \,\otimes\, g,\, S^{\,-1}_{C_{1}}\,f_{\,i} \,\otimes\, S^{\,-1}_{C_{2}}\,g_{\,j} \,|\, a_{\,2} \otimes b_{\,2}, \,\cdots,\, a_{\,n} \otimes b_{\,n}\,\right>\,\left(\,C_{1}\,f_{\,i} \,\otimes\, C_{2}\,g_{\,j}\,\right)\\
&=\sum\limits^{\infty}_{i,\,j \,=\, 1}\,\left <\,f \,\otimes\, g,\, \left(\,S^{\,-1}_{C_{1}} \,\otimes\, S^{\,-1}_{C_{2}}\,\right)\,\left(\,f_{\,i} \,\otimes\, g_{\,j}\,\right) \,|\, a_{\,2} \otimes b_{\,2}, \,\cdots,\, a_{\,n} \otimes b_{\,n}\,\right>\,\left(\,C_{1}\,f_{\,i} \,\otimes\, C_{2}\,g_{\,j}\,\right)\\
 &\,=\, \sum\limits^{\infty}_{i,\,j \,=\, 1}\,\left <\,f \,\otimes\, g,\, S^{\,-1}_{C_{1} \,\otimes\, C_{2}}\,\left(\,f_{\,i} \,\otimes\, g_{\,j}\,\right) \,|\, a_{\,2} \otimes b_{\,2}, \,\cdots,\, a_{\,n} \otimes b_{\,n}\,\right >\,\left(\,C_{1} \,\otimes\, C_{2}\,\right)\,\left(\,f_{\,i} \,\otimes\, g_{\,j}\,\right).
\end{align*}
Since \,$\left\{\,\left <\,f,\, S^{\,-1}_{C_{1}}\,f_{\,i} \,|\, a_{\,2},\, \cdots,\, a_{\,n}\,\,\right >_{1}\,\right \}^{\infty}_{i \,=\, 1},\, \left\{\,\left <\,g,\, S^{\,-1}_{C_{2}}\,g_{\,j} \,|\, b_{\,2},\, \cdots,\, b_{\,n}\,\,\right >_{2}\,\right \}^{\infty}_{j \,=\, 1} \,\in\, l^{\,2} (\,\mathbb{N}\,)$, \,$\left\{\,f_{\,i}\,\right\}^{\infty}_{i \,=\, 1}$\, and \,$\left\{\,g_{\,j}\,\right\}^{\infty}_{j \,=1\, }$\, are \,$C_{1}$-controlled and \,$C_{2}$-controlled Bessel sequence associated to \,$\left(\,a_{\,2},\, \cdots,\, a_{\,n}\,\right)$\, and \,$\left(\,b_{\,2},\, \cdots,\, b_{\,n}\,\right)$\, in \,$H$\, and \,$K$, respectively, the above series converges unconditionally.\,On the other hand, for each \,$f \,\otimes\, g \,\in\, H_{F} \,\otimes\, K_{G}$, 
\begin{align*}
&f \,\otimes\, g \,=\, S^{\,-1}_{C_{1}}\, S_{C_{1}}\,f \,\otimes\, S^{\,-1}_{C_{2}}\, S_{C_{2}}\,g \\
&=\,S^{\,-1}_{C_{1}}\,\left(\,\sum\limits^{\infty}_{i \,=\, 1}\,\left <\,f,\, f_{\,i} \,|\, a_{\,2},\, \cdots,\, a_{\,n}\,\,\right >_{1}\,C_{1}\,f_{\,i}\,\right) \,\otimes\, S^{\,-1}_{C_{2}}\,\left(\,\sum\limits^{\infty}_{j \,=\, 1}\,\left <\,g,\, g_{\,j} \,|\, b_{\,2},\, \cdots,\, b_{\,n}\,\,\right >_{2}\,C_{2}\,g_{\,j}\,\right)\\
 &\,=\, \sum\limits^{\infty}_{i,\,j \,=\, 1}\,\left <\,f \,\otimes\, g,\, f_{\,i} \,\otimes\, g_{\,j} \,|\, a_{\,2} \otimes b_{\,2}, \,\cdots,\, a_{\,n} \otimes b_{\,n}\,\right >\,S^{\,-1}_{C_{1} \,\otimes\, C_{2}}\,\left(\,C_{1} \,\otimes\, C_{2}\,\right)\,\left(\,f_{\,i} \,\otimes\, g_{\,j}\,\right).  
\end{align*}
This completes the proof.    
\end{proof}

\begin{corollary}
Let \,$\left\{\,f_{\,i}\,\right\}^{\infty}_{i \,=1\, }$\, and \,$\left\{\,g_{\,j}\,\right\}^{\infty}_{j \,=1\, }$\, be \,$C_{1}$-controlled and \,$C_{2}$-controlled tight frames associated to \,$\left(\,a_{\,2},\, \cdots,\, a_{\,n}\,\right)$\, and \,$\left(\,b_{\,2},\, \cdots,\, b_{\,n}\,\right)$\, for \,$H$\, and \,$K$\, with bounds \,$A_{1}$\, and \,$A_{2}$, respectively.\,Then for each \,$f \,\otimes\, g \,\in\, H_{F} \,\otimes\, K_{G}$, we have
\[ f  \,\otimes\, g \,= \dfrac{1}{A_{1}\,A_{2}}\sum\limits^{\infty}_{i,\,j \,=\, 1}\left<\,f \,\otimes\, g,\, f_{\,i} \,\otimes\, g_{\,j} \,|\, a_{\,2} \otimes b_{\,2}, \,\cdots,\, a_{\,n} \otimes b_{\,n}\,\right >\left(\,C_{1} \,\otimes\, C_{2}\,\right)\left(\,f_{\,i} \,\otimes\, g_{\,j}\,\right).\]
\end{corollary}

In the next theorem, we establish that an image of controlled frame associated to \,$\left(\,a_{\,2},\, \cdots,\, a_{\,n}\,\right)$\, under a bounded linear operator becomes a controlled frame associated to \,$\left(\,a_{\,2},\, \cdots,\, a_{\,n}\,\right)$\, if and only if the bounded linear operator have to be invertible.

\begin{theorem}
Let \,$\{\,f_{\,i}\,\}_{i \,=\,1}^{\infty}$\, be a \,$C_{1}$-controlled frame associated to \,$\left(\,a_{\,2},\, \cdots,\, a_{\,n}\,\right)$\, for \,$H$\, with bounds \,$A,\,B$\, and \,$\{\,g_{\,j}\,\}_{j \,=\,1}^{\infty}$\, be a \,$C_{2}$-controlled frame associated to \,$(\,b_{\,2},\, \cdots$, \,$b_{\,n}\,)$\, for \,$K$\, with bounds \,$C,\,D$\, having their corresponding frame operators \,$S_{C_{1}}$\, and \,$S_{C_{2}}$, respectively.\,Suppose \,$C_{1}$\, and \,$C_{2}$\, commutes with \,$U_{1}$\, and \,$U_{2}$, respectively.\,Then \,$ \left\{\,\Delta_{i\,j} \,=\, \left(\,U_{1} \,\otimes\, U_{2}\,\right)\,\left(\,f_{\,i} \,\otimes\, g_{\,j}\,\right)\,\right\}^{\infty}_{i,\,j \,=\, 1}$\, is a \,$\left(\,C_{1} \,\otimes\, C_{2}\,\right)$-controlled frame associated to \,$(\,a_{\,2} \,\otimes\, b_{\,2},\, \,\cdots$,\, $a_{\,n} \,\otimes\, b_{\,n}\,)$\, for \,$H \,\otimes\, K$\, if and only if \,$U_{1} \,\otimes\, U_{2}$\, is an invertible operator on \,$H_{F} \,\otimes\, K_{G}$.
\end{theorem}

\begin{proof}
First we suppose that \,$U_{1} \,\otimes\, U_{2}$\, is an invertible on \,$H_{F} \,\otimes\, K_{G}$.\;Then by Theorem \ref{th1.1}, \,$U_{1}$\, and \,$U_{2}$\, are invertible on \,$H_{F}$\, and \,$K_{G}$, respectively.\;For each \,$f \,\in\, H_{F}$\, and \,$g \,\in\, K_{G}$, we have
\begin{equation}\label{eq1.05} 
\left\|\,f,\, a_{\,2},\, \cdots,\, a_{\,n}\,\right\|_{1} \,\leq\, \left\|\,U^{\,-\,1}_{1}\,\right\|\,\left\|\,U^{\,\ast}_{1}\,(\,f\,),\, a_{\,2},\, \cdots,\, a_{\,n}\,\right\|_{1},\;\text{and}
\end{equation}
\begin{equation}\label{eq1.5}
\left\|\,g,\, b_{\,2},\, \cdots,\, b_{\,n}\,\right\|_{2} \,\leq\, \left\|\,U^{\,-\, 1}_{2}\,\right\|\,\left\|\,U^{\,\ast}_{2}\,(\,g\,),\, b_{\,2},\, \cdots,\, b_{\,n}\,\right\|_{2}.\hspace{.7cm}
\end{equation}
Now, for each  \,$f \,\otimes\, g \,\in\, H_{F} \,\otimes\, K_{G}$, we have
\begin{align*}
&\sum\limits_{i,\, j \,=\, 1}^{\,\infty}\,\left<\,f \,\otimes\, g,\, \left(\,U_{1} \,\otimes\, U_{2}\,\right)\,\left(\,f_{\,i} \,\otimes\, g_{\,j}\,\right) \,|\, a_{\,2} \,\otimes\, b_{\,2},\, \cdots,\, a_{\,n} \,\otimes\, b_{\,n}\,\right>\,\times\\
&\hspace{1.2cm}\left<\,\left(\,C_{1} \,\otimes\, C_{2}\,\right)\,\left(\,U_{1} \,\otimes\, U_{2}\,\right)\,\left(\,f_{\,i} \,\otimes\, g_{\,j}\,\right),\, f \,\otimes\, g \,|\, a_{\,2} \,\otimes\, b_{\,2},\, \cdots,\, a_{\,n} \,\otimes\, b_{\,n}\,\right>\\
&=\,\sum\limits_{i,\, j \,=\, 1}^{\,\infty}\,\left<\,f \,\otimes\, g,\, \left(\,U_{1}\,f_{\,i} \,\otimes\, U_{2}\,g_{\,j}\,\right) \,|\, a_{\,2} \,\otimes\, b_{\,2},\, \cdots,\, a_{\,n} \,\otimes\, b_{\,n}\,\right>\,\times\\
&\hspace{1.2cm}\left<\,\left(\,C_{1}\,U_{1}\,f_{\,i} \,\otimes\, C_{2}\,U_{2}\,g_{\,j}\,\right),\, f \,\otimes\, g \,|\, a_{\,2} \,\otimes\, b_{\,2},\, \cdots,\, a_{\,n} \,\otimes\, b_{\,n}\,\right>\\
&=\,\sum\limits_{i,\, j \,=\, 1}^{\,\infty}\,\left<\,f,\, U_{1}\,f_{\,i}\,|\,a_{\,2},\, \cdots,\, a_{\,n}\,\right>_{1}\,\left<\,g,\, U_{2}\,g_{\,j} \,|\, b_{\,2},\, \cdots,\, b_{\,n}\,\right>_{2}\,\times\\
&\hspace{1.5cm}\left<\,C_{1}\,U_{1}\,f_{\,i},\, f \,|\,a_{\,2},\, \cdots,\, a_{\,n}\,\right>_{1}\,\left<\,C_{2}\,U_{2}\,g_{\,j},\, g \,|\, b_{\,2},\, \cdots,\, b_{\,n}\,\right>_{2}\\
&=\,\sum\limits_{i,\, j \,=\, 1}^{\,\infty}\,\left<\,U_{1}^{\,\ast}\,f,\, f_{\,i}\,|\,a_{\,2},\, \cdots,\, a_{\,n}\,\right>_{1}\,\left<\,U_{2}^{\,\ast}\,g,\, g_{\,j} \,|\, b_{\,2},\, \cdots,\, b_{\,n}\,\right>_{2}\,\times\\
&\hspace{1.5cm}\left<\,C_{1}\,f_{\,i},\, U_{1}^{\,\ast}\,f \,|\,a_{\,2},\, \cdots,\, a_{\,n}\,\right>_{1}\,\left<\,C_{2}\,g_{\,j},\, U_{2}^{\,\ast}\,g \,|\, b_{\,2},\, \cdots,\, b_{\,n}\,\right>_{2}\\
&=\,\sum\limits_{\,i \,=\, 1}^{\,\infty}\,\left<\,U_{1}^{\,\ast}\,f,\, f_{\,i} \,|\, a_{\,2},\, \cdots,\, a_{\,n}\,\right>_{1}\,\left<\,C_{1}\,f_{\,i},\, U_{1}^{\,\ast}\,f \,|\,a_{\,2},\, \cdots,\, a_{\,n}\,\right>_{1}\,\times\\
&\hspace{1.5cm}\sum\limits_{\,j \,=\, 1}^{\,\infty}\,\left<\,U_{2}^{\,\ast}\,g,\, g_{\,j} \,|\, b_{\,2},\, \cdots,\, b_{\,n}\,\right>_{2}\,\left<\,C_{2}\,g_{\,j},\, U_{2}^{\,\ast}\,g \,|\, b_{\,2},\, \cdots,\, b_{\,n}\,\right>_{2}\\
&\leq\,B\,D\,\left\|\,U_{1}^{\,\ast}\,f,\, a_{\,2},\, \cdots,\, a_{\,n}\,\right\|_{1}^{\,2}\,\left\|\,U_{2}^{\,\ast}\,g,\, b_{\,2},\, \cdots,\, b_{\,n}\,\right\|_{2}^{\,2}\\
&\leq\, B\,D\,\left\|\,U_{1}^{\,\ast}\,\right\|^{\,2}\,\left\|\,U_{2}^{\,\ast}\,\right\|^{\,2}\,\left\|\,f,\, a_{\,2},\, \cdots,\, a_{\,n}\,\right\|_{1}^{\,2}\,\left\|\,g \,,\, b_{\,2},\, \cdots,\, b_{\,n}\,\right\|_{2}^{\,2}\\
&=\, B\,D\,\left\|\,U_{1} \,\otimes\, U_{2}\,\right\|^{\,2}\,\left\|\,f \,\otimes\, g,\, a_{\,2} \,\otimes\, b_{\,2} \,\cdots,\, a_{\,n} \,\otimes\, b_{\,n}\,\right\|^{\,2}. 
\end{align*} 
On the other hand,
\begin{align*}
&\sum\limits_{i,\, j \,=\, 1}^{\,\infty}\,\left<\,f \,\otimes\, g,\, \left(\,U_{1} \,\otimes\, U_{2}\,\right)\,\left(\,f_{\,i} \,\otimes\, g_{\,j}\,\right) \,|\, a_{\,2} \,\otimes\, b_{\,2},\, \cdots,\, a_{\,n} \,\otimes\, b_{\,n}\,\right>\,\times\\
&\hspace{1.2cm}\left<\,\left(\,C_{1} \,\otimes\, C_{2}\,\right)\,\left(\,U_{1} \,\otimes\, U_{2}\,\right)\,\left(\,f_{\,i} \,\otimes\, g_{\,j}\,\right),\, f \,\otimes\, g \,|\, a_{\,2} \,\otimes\, b_{\,2},\, \cdots,\, a_{\,n} \,\otimes\, b_{\,n}\,\right>\\
&\geq\, A\,C\,\left\|\,U_{1}^{\,\ast}\,f,\, a_{\,2},\, \cdots,\, a_{\,n}\,\right\|_{1}^{\,2}\,\left\|\,U_{2}^{\,\ast}\,f,\, b_{\,2},\, \cdots,\, b_{\,n}\,\right\|_{2}^{\,2}\\
&\geq\, \dfrac{A\,C\left\|\,f,\, a_{\,2},\, \cdots,\, a_{\,n}\,\right\|_{1}^{\,2}\,\left\|\,g,\, b_{\,2},\, \cdots,\, b_{\,n}\,\right\|_{2}^{\,2}}{\left\|\,U^{\,-\,1}_{1}\,\right\|^{\,2}\,\left\|\,U^{\,-\,1}_{2}\,\right\|^{\,2}}\;[\;\text{by (\ref{eq1.05}) and (\ref{eq1.5})}\;]\\
&=\, \dfrac{A\,C}{\left\|\,\left(\,U_{1} \,\otimes\, U_{2}\,\right)^{\,-\, 1}\,\right\|^{\,2}}\,\left\|\,f \,\otimes\, g,\, a_{\,2} \,\otimes\, b_{\,2},\,\cdots,\, a_{\,n} \,\otimes\, b_{\,n}\,\right\|^{\,2}.
\end{align*}
Therefore, the sequence \,$\left\{\,\Delta_{i\,j}\,\right\}_{i,\,j \,=\, 1}^{\,\infty}$\, is a \,$\left(\,C_{1} \,\otimes\, C_{2}\,\right)$-controlled frame associated to \,$(\,a_{\,2} \,\otimes\, b_{\,2},\, \,\cdots,\, a_{\,n} \,\otimes\, b_{\,n}\,)$\, for \,$H \,\otimes\, K$.\\

Conversely, suppose that \,$\left\{\,\Delta_{i\,j}\,\right\}_{i,\,j \,=\, 1}^{\,\infty}$\, is a \,$\left(\,C_{1} \,\otimes\, C_{2}\,\right)$-controlled frame associated to \,$(\,a_{\,2} \,\otimes\, b_{\,2},\, \,\cdots$, \,$a_{\,n} \,\otimes\, b_{\,n}\,)$\, for \,$H \,\otimes\, K$.\;Now, for each \,$f \,\otimes\, g \,\in\, H_{F} \,\otimes\, K_{G}$,
\begin{align*}
&\sum\limits_{i,\, j \,=\, 1}^{\,\infty}\,\left<\,f \,\otimes\, g,\,  \Delta_{i\,j} \,|\, a_{\,2} \,\otimes\, b_{\,2}, \,\cdots,\, a_{\,n} \,\otimes\, b_{\,n}\,\right>\,\left(\,C_{1} \,\otimes\, C_{2}\,\right)\,\Delta_{i\,j}\\
&=\, \sum\limits_{i,\, j \,=\, 1}^{\,\infty}\left<\,f \otimes g,\, U_{1}\,f_{\,i} \,\otimes\, U_{2}\,g_{\,j} \,|\,a_{\,2} \,\otimes\, b_{\,2}, \,\cdots,\, a_{\,n} \otimes b_{\,n}\,\right>\,\left(\,C_{1}\,U_{1}\,f_{\,i} \,\otimes\, C_{2}\,U_{2}\,g_{\,j}\,\right)\\
&= \left(\sum\limits_{\,i \,=\, 1}^{\,\infty}\left<\,U_{1}^{\,\ast}f,\, f_{i}\,|\,a_{2},\, \cdots,\, a_{n}\,\right>_{1}\,U_{1}\,C_{1}\,f_{i}\,\right) \,\otimes\, \left(\sum\limits_{\,j \,=\, 1}^{\,\infty}\left<\,U_{2}^{\,\ast}\,g,\, g_{j}\,|\,b_{2},\, \cdots,\, b_{n}\,\right>_{2}\,U_{2}\,C_{2}\,g_{j}\right)\\
&=\, U_{1}\,S_{C_{1}}\,U_{1}^{\,\ast}\,f \,\otimes\,  U_{2}\,S_{C_{2}}\,U_{2}^{\,\ast}\,g = \left(\,U_{1} \,\otimes\, U_{2}\,\right)\,\left(\,S_{C_{1}} \,\otimes\, S_{C_{2}}\,\right)\,\left(\,U^{\,\ast}_{1} \,\otimes\, U^{\,\ast}_{2}\right)\,(\,f \,\otimes\, g\,)\\
&=\, \left(\,U_{1} \,\otimes\, U_{2}\,\right)\,S_{C_{1} \,\otimes\, C_{2}}\,\left(\,U_{1} \,\otimes\, U_{2}\,\right)^{\,\ast}\,(\,f \,\otimes\, g\,)
\end{align*} 
Hence, the frame operator for \,$\left\{\,\Delta_{i\,j}\,\right\}_{i,\,j \,=\, 1}^{\,\infty}$\, is \,$\left(\,U_{1} \,\otimes\, U_{2}\,\right)\,S_{C_{1} \,\otimes\, C_{2}}\,\left(\,U_{1} \,\otimes\, U_{2}\,\right)^{\,\ast}$\, and therefore it is invertible.\;Also, we know that \,$S_{C_{1} \,\otimes\, C_{2}}$\, is invertible and hence \,$U_{1} \,\otimes\, U_{2}$\, is invertible on \,$H_{F} \,\otimes\, K_{G}$.\,This completes the proof.
\end{proof}

We now present the concept of a dual controlled frame in \,$H \,\otimes\, K$.

\begin{definition}
Let \,$\left\{\,f_{\,i} \otimes g_{\,j}\,\right\}^{\,\infty}_{i,\,j \,=\, 1}$\, be a \,$\left(\,C_{1} \,\otimes\, C_{2}\,\right)$-controlled frame associated to \,$(\,a_{\,2} \otimes b_{\,2},\, \cdots$, \,$a_{\,n} \otimes b_{\,n}\,)$\, for \,$H \otimes K$.\,Then a \,$\left(\,C_{1} \,\otimes\, C_{2}\,\right)$-controlled frame \,$\left\{\,e_{\,i} \otimes h_{\,j}\,\right\}^{\,\infty}_{i,\,j \,=\, 1}$\, associated to \,$\left(\,a_{\,2} \,\otimes\, b_{\,2},\, \,\cdots,\, a_{\,n} \,\otimes\, b_{\,n}\,\right)$\, for \,$H \,\otimes\, K$\, satisfying 
\begin{align}
&f \,\otimes\, g \nonumber\\
&\,=\, \sum\limits_{i,\, j \,=\, 1}^{\,\infty}\,\left<\,f \,\otimes\, g,\, e_{\,i} \,\otimes\, h_{\,j} \,|\, a_{\,2} \,\otimes\, b_{\,2}, \,\cdots,\, a_{\,n} \,\otimes\, b_{\,n}\,\right>\,\left(\,C_{1} \,\otimes\, C_{2}\,\right)\,\left(\,f_{\,i} \,\otimes\, g_{\,j}\,\right),\label{eq2.1}
\end{align}
for all \,$f \,\otimes\, g \,\in\, H \,\otimes\, K$, is called a dual \,$\left(\,C_{1} \,\otimes\, C_{2}\,\right)$-controlled frame associated to \,$(\,a_{\,2} \,\otimes\, b_{\,2},\, \cdots$, \,$a_{\,n} \,\otimes\, b_{\,n}\,)$\, for \,$H \otimes K$\, of \,$\left\{\,f_{\,i} \,\otimes\, g_{\,j}\,\right\}^{\infty}_{i,\,j = 1}$.     
\end{definition}

\begin{theorem}\label{th3.2}
Let \,$\{\,f_{\,i}\,\}_{i \,=\,1}^{\infty} \,,\, \left\{\,e_{\,i}\,\right\}^{\,\infty}_{i \,=\, 1}$\, be a pair of dual \,$C_{1}$-controlled frames associated to \,$\left(\,a_{\,2},\, \cdots,\, a_{\,n}\,\right)$\, for \,$H$\, and \,$\{\,g_{\,j}\,\}_{j \,=\,1}^{\infty} \,,\, \left\{\,h_{\,j}\,\right\}^{\,\infty}_{j \,=\, 1}$\, be a pair of dual \,$C_{2}$-controlled frames associated to \,$\left(\,b_{\,2},\, \cdots,\, b_{\,n}\,\right)$\, for \,$K$.\;Then \,$\left\{\,e_{\,i} \,\otimes\, h_{\,j}\,\right\}^{\,\infty}_{i,\,j \,=\, 1}$\, is a dual \,$\left(\,C_{1} \,\otimes\, C_{2}\,\right)$-controlled frame associated to \,$\left(\,a_{\,2} \,\otimes\, b_{\,2},\, \,\cdots,\, a_{\,n} \,\otimes\, b_{\,n}\,\right)$\, for \,$H \otimes K$\, of \,$\left\{\,f_{\,i} \,\otimes\, g_{\,j}\,\right\}^{\,\infty}_{i,\,j \,=\, 1}$.     
\end{theorem}

\begin{proof}
By Theorem \ref{th2.1}, \,$\left\{\,f_{\,i} \,\otimes\, g_{\,j}\,\right\}^{\,\infty}_{i,\,j \,=\, 1}$, \,$\left\{\,e_{\,i} \,\otimes\, h_{\,j}\,\right\}^{\,\infty}_{i,\,j \,=\, 1}$\, are \,$\left(\,C_{1} \,\otimes\, C_{2}\,\right)$-controlled frames associated to \,$\left(\,a_{2} \,\otimes\, b_{2},\, \,\cdots,\, a_{n} \,\otimes\, b_{n}\,\right)$\, for \,$H \,\otimes\, K$.\;Since \,$\left\{\,e_{\,i}\,\right\}^{\,\infty}_{i \,=\, 1}$\, and \,$\left\{\,h_{\,j}\,\right\}^{\,\infty}_{j \,=\, 1}$\, are dual \,$C_{1}$-controlled and \,$C_{2}$-controlled frames associated to \,$\left(a_{2},\, \cdots,\, a_{n}\right)$\, and \,$\left(b_{2},\, \cdots,\, b_{n}\right)$\, of \,$\{\,f_{\,i}\,\}_{i \,=\,1}^{\infty}$\, and \,$\{\,g_{\,j}\,\}_{j \,=\,1}^{\infty}$, respectively, for all \,$f \,\in\, H_{F}$, \,$g \,\in\, K_{G}$, 
\[f \,=\, \sum\limits^{\,\infty}_{i \,=\, 1}\left<\,f,\, e_{\,i} \,|\, a_{\,2},\, \cdots,\, a_{\,n}\,\right>_{1}\,C_{1}\,f_{\,i},\, \;\text{and}\; g \,=\, \sum\limits^{\,\infty}_{j \,=\, 1}\left<\,g,\, h_{j} \,|\, b_{\,2},\, \cdots,\, b_{\,n}\,\right>_{2}\,C_{2}\,g_{j}.\] 
Then, for all \,$f \,\otimes\, g \,\in\, H_{F} \,\otimes\, K_{G}$, we have
\begin{align*}
&f \,\otimes\, g \\
&\,=\, \left(\sum\limits^{\,\infty}_{i \,=\, 1}\left<\,f,\, e_{\,i} \,|\, a_{\,2},\, \cdots,\, a_{\,n}\,\right>_{1}\,C_{1}\,f_{\,i}\right) \otimes \left(\sum\limits^{\,\infty}_{j \,=\, 1}\left<\,g,\, h_{j} \,|\, b_{\,2},\, \cdots,\, b_{\,n}\,\right>_{2}\,C_{2}\,g_{j}\right)\\
&=\, \sum\limits_{i,\, j \,=\, 1}^{\,\infty}\,\left<\,f,\, e_{\,i} \,|\, a_{\,2},\, \cdots,\, a_{\,n}\,\right>_{1}\, \left<\,g,\, h_{\,j} \,|\, b_{\,2},\, \cdots,\, b_{\,n}\,\right>_{2}\,\left(\,C_{1}\,f_{\,i} \,\otimes\, C_{2}\,g_{\,j}\,\right)\\
&=\, \sum\limits_{i,\, j \,=\, 1}^{\,\infty}\,\left<\,f \,\otimes\, g,\, e_{\,i} \,\otimes\, h_{\,j} \,|\, a_{\,2} \,\otimes\, b_{\,2}, \,\cdots,\, a_{\,n} \,\otimes\, b_{\,n}\,\right>\,\left(\,C_{1} \,\otimes\, C_{2}\,\right)\,\left(\,f_{\,i} \,\otimes\, g_{\,j}\,\right).
\end{align*}
This completes the proof.    
\end{proof}

\begin{theorem}
Let \,$\{\,f_{\,i}\,\}_{i \,=\,1}^{\infty},\, \left\{\,e_{\,i}\,\right\}^{\,\infty}_{i \,=\, 1}$\, be a pair of dual \,$C_{1}$-controlled frames associated to \,$\left(a_{\,2},\, \cdots,\, a_{\,n}\right)$\, for \,$H$\, and \,$\{\,g_{\,j}\,\}_{j \,=\,1}^{\infty},\, \left\{\,h_{\,j}\,\right\}^{\,\infty}_{j \,=\, 1}$\, be a pair of dual \,$C_{2}$-controlled frames associated to \,$\left(b_{\,2},\, \cdots,\, b_{\,n}\right)$\, for \,$K$.\,Suppose \,$U \,\in\, \mathcal{B}\,(\,H_{F}\,),\, \,V \,\in\, \mathcal{B}\,(\,K_{G}\,)$\, are unitary operators, \,$C_{1}$\, and \,$C_{2}$\, commutes with \,$U$\, and \,$V$, respectively.\,Then \,$\Lambda \,=\, \left\{\,\left(\,U \,\otimes\, V\,\right)\,\left(\,f_{\,i} \,\otimes\, g_{\,j}\,\right)\,\right\}^{\,\infty}_{i,\,j = 1}$\, and \,$\Gamma \,=\, \left\{\,\left(\,U \,\otimes\, V\,\right)\,\left(\,e_{\,i} \,\otimes\, h_{\,j}\,\right)\,\right\}^{\,\infty}_{i,\,j \,=\, 1}$\, also form a pair of dual \,$\left(\,C_{1} \,\otimes\, C_{2}\,\right)$-controlled frames associated to \,$(\,a_{\,2} \,\otimes\, b_{\,2},\, \,\cdots,\, a_{\,n} \,\otimes\, b_{\,n}\,)$\, for \,$H \,\otimes\, K$.   
\end{theorem}

\begin{proof}
By Theorem \ref{th3.2}, the sequences \,$\left\{\,f_{\,i} \,\otimes\, g_{\,j}\,\right\}^{\,\infty}_{i,\,j \,=\, 1}$\, and \,$\left\{\,e_{\,i} \,\otimes\, h_{\,j}\,\right\}^{\,\infty}_{i,\,j \,=\, 1}$\, form a pair of dual \,$\left(\,C_{1} \,\otimes\, C_{2}\,\right)$-controlled frames associated to \,$\left(\,a_{\,2} \,\otimes\, b_{\,2},\, \,\cdots,\, a_{\,n} \,\otimes\, b_{\,n}\,\right)$\, for \,$H \,\otimes\, K$.\;Now, for each \,$f \,\otimes\, g \,\in\, H_{F} \,\otimes\, K_{G}$, we have
\begin{align*}
&\sum\limits_{i,\, j \,=\, 1}^{\,\infty}\,\left<\,f \,\otimes\, g,\, \left(\,U \,\otimes\, V\,\right)\,\left(\,f_{\,i} \,\otimes\, g_{\,j}\,\right) \,|\, a_{\,2} \,\otimes\, b_{\,2},\, \cdots,\, a_{\,n} \,\otimes\, b_{\,n}\,\right>\,\times\\
&\hspace{1cm}\left<\,\left(\,C_{1} \,\otimes\, C_{2}\,\right)\,\left(\,U \,\otimes\, V\,\right)\,\left(\,f_{\,i} \,\otimes\, g_{\,j}\,\right),\, f \,\otimes\, g \,|\, a_{\,2} \,\otimes\, b_{\,2},\, \cdots,\, a_{\,n} \,\otimes\, b_{\,n}\,\right>\\
&=\,\sum\limits_{i,\, j \,=\, 1}^{\,\infty}\,\left<\,f \,\otimes\, g,\, \left(\,U\,f_{\,i} \,\otimes\, V\,g_{\,j}\,\right) \,|\, a_{\,2} \,\otimes\, b_{\,2},\, \cdots,\, a_{\,n} \,\otimes\, b_{\,n}\,\right>\,\times\\
&\hspace{1cm}\left<\,\left(\,C_{1}\,U\,f_{\,i} \,\otimes\, C_{2}\,V\,g_{\,j}\,\right),\, f \,\otimes\, g \,|\, a_{\,2} \,\otimes\, b_{\,2},\, \cdots,\, a_{\,n} \,\otimes\, b_{\,n}\,\right>\\
&= \sum\limits_{\,i \,=\, 1}^{\,\infty}\,\left<\,f,\, U\,f_{i}\,|\,a_{2},\, \cdots,\, a_{n}\,\right>_{1}\,\left<\,C_{1}\,U\,f_{\,i},\, f \,|\,a_{2},\, \cdots,\, a_{n}\,\right>_{1}\,\times\\
&\hspace{1cm}\sum\limits_{\,j \,=\, 1}^{\,\infty}\,\left<\,g,\, V\,g_{\,j}\,|\, b_{2},\, \cdots,\, b_{n}\,\right>_{2}\,\left<\,C_{2}\,V\,g_{\,j},\, g \,|\, b_{2},\, \cdots,\, b_{n}\,\right>_{2}\\
&=\,\sum\limits_{\,i \,=\, 1}^{\,\infty}\,\left<\,U^{\,\ast}\,f,\,  f_{i}\,|\,a_{2},\, \cdots,\, a_{n}\,\right>_{1}\,\left<\,C_{1}\,f_{\,i},\, U^{\,\ast}\,f \,|\,a_{2},\, \cdots,\, a_{n}\,\right>_{1}\,\times\\
&\hspace{1cm}\sum\limits_{\,j \,=\, 1}^{\,\infty}\,\left<\,V^{\,\ast}\,g,\, g_{\,j}\,|\, b_{2},\, \cdots,\, b_{n}\,\right>_{2}\,\left<\,C_{2}\,g_{\,j},\, V^{\,\ast}\,g \,|\, b_{2},\, \cdots,\, b_{n}\,\right>_{2}.\\ 
\end{align*}
Since \,$\{\,f_{\,i}\,\}_{i \,=\,1}^{\infty}$\, is a \,$C_{1}$-controlled frame associated to \,$\left(\,a_{\,2},\, \cdots,\, a_{\,n}\,\right)$\, for \,$H$\, and \,$\{\,g_{\,j}\,\}_{j \,=\,1}^{\infty}$\, is a \,$C_{2}$-controlled frame associated to \,$\left(\,b_{\,2},\, \cdots,\, b_{\,n}\,\right)$\, for \,$K$, the above calculation shows that the sequence \,$\Lambda$\, is a \,$\left(\,C_{1} \,\otimes\, C_{2}\,\right)$-controlled frame associated to \,$\left(\,a_{\,2} \,\otimes\, b_{\,2},\, \,\cdots,\, a_{\,n} \,\otimes\, b_{\,n}\,\right)$\, for \,$H \,\otimes\, K$.\;Similarly, it can be shown that \,$\Gamma$\, is a \,$\left(\,C_{1} \,\otimes\, C_{2}\,\right)$-controlled frame associated to \,$\left(\,a_{\,2} \,\otimes\, b_{\,2},\, \,\cdots,\, a_{\,n} \,\otimes\, b_{\,n}\,\right)$\, for \,$H \,\otimes\, K$.  
Furthermore, for each \,$f \,\otimes\, g \,\in\, H_{F} \,\otimes\, K_{G}$, we have
\begin{align*}
&\sum\limits_{i,\, j = 1}^{\infty}\left<\,f \otimes f,\, \left(U \otimes V\right)\,\left(e_{i} \otimes h_{j}\right)\,|\,a_{2} \otimes b_{2}, \,\cdots,\, a_{n} \otimes b_{n}\,\right>\left(\,C_{1} \,\otimes\, C_{2}\,\right)\left(U \otimes V\right)\left(f_{i} \otimes g_{j}\right)\\
&=\, \sum\limits_{i,\, j \,=\, 1}^{\,\infty}\left<\,f \,\otimes\, g,\, \left(\,U\,e_{i} \,\otimes\, V\,h_{j}\,\right) \,|\, a_{2} \,\otimes\, b_{2}, \,\cdots,\, a_{n} \,\otimes\, b_{n}\,\right>\,\left(\,C_{1}\,U\,f_{i} \,\otimes\, C_{2}\,V\,g_{j}\,\right)\\
&=\, \sum\limits_{i,\, j \,=\, 1}^{\,\infty}\left<\,f \,\otimes\, g,\, \left(\,U\,e_{i} \,\otimes\, V\,h_{j}\,\right) \,|\, a_{2} \,\otimes\, b_{2}, \,\cdots,\, a_{n} \,\otimes\, b_{n}\,\right>\,\left(\,U\,C_{1}\,f_{i} \,\otimes\, V\,C_{2}\,g_{j}\,\right)\\
&= U\sum\limits_{\,i = 1}^{\infty}\left<\,U^{\,\ast}\,f,\, e_{\,i}\,|\,a_{\,2},\, \cdots,\, a_{\,n}\,\right>_{1}\,C_{1}\,f_{\,i} \,\otimes\, V\sum\limits_{\,j = 1}^{\,\infty}\left<\,V^{\,\ast}\,g,\, h_{\,j}\,|\,b_{\,2},\, \cdots,\, b_{\,n}\,\right>_{2}\,C_{2}\,g_{\,j}\\
& \,=\, U\,U^{\,\ast}\,(\,f\,) \,\otimes\, V\,V^{\,\ast}\,(\,g\,) \,=\, f \,\otimes\, g.
\end{align*}
Hence, \,$\Lambda$\, and \,$\Gamma$\, form a pair of dual \,$\left(\,C_{1} \,\otimes\, C_{2}\,\right)$-controlled frames associated to \,$\left(\,a_{\,2} \,\otimes\, b_{\,2},\, \,\cdots,\, a_{\,n} \,\otimes\, b_{\,n}\,\right)$\, for \,$H \,\otimes\, K$.
\end{proof}

Now, we end this section by considering controlled frame in the direct sum of \,$n$-Hilbert spaces \,$H$\, and \,$K$.\,The direct sum of \,$n$-Hilbert spaces \,$H$\, and \,$K$\, is denoted by \,$H \,\oplus\, K$\, and defined to be a \,$n$-Hilbert space associated with the \,$n$-inner product 
\[\left<\,f_{\,1} \,\oplus\, g_{\,1},\, f_{\,2} \,\oplus\, g_{\,2} \,|\, f_{\,3} \,\oplus\, g_{\,3},\, \,\cdots,\, f_{\,n} \,\oplus\, g_{\,n}\,\right>\]
\begin{equation}\label{eqnn1}
 \,=\, \left<\,f_{\,1},\, f_{\,2} \,|\, f_{\,3},\, \,\cdots,\, f_{\,n}\,\right>_{1} \,+\, \left<\,g_{\,1},\, g_{\,2} \,|\, g_{\,3},\, \,\cdots,\, g_{\,n}\,\right>_{2},
\end{equation}
for all \,$f_{\,1},\, f_{\,2},\, f_{\,3},\, \,\cdots,\, f_{\,n} \,\in\, H$\, and \,$g_{\,1},\, g_{\,2},\, g_{\,3},\, \,\cdots,\, g_{\,n} \,\in\, K$.\\
The \,$n$-norm on \,$H \,\otimes\, K$\, is defined by 
\[\left\|\,f_{\,1} \,\oplus\, g_{\,1},\, f_{\,2} \,\oplus\, g_{\,2},\, \,\cdots,\,\, f_{\,n} \,\oplus\, g_{\,n}\,\right\|\]
\begin{equation}\label{eqnn1.1}
\hspace{.6cm} =\,\left\|\,f_{\,1},\, f_{\,2},\, \cdots,\, f_{\,n}\,\right\|_{1} \,+\, \left\|\,g_{\,1},\, g_{\,2},\, \cdots,\, g_{\,n}\,\right\|_{2},
\end{equation}
for all \,$f_{\,1},\, f_{\,2},\, \,\cdots,\, f_{\,n} \,\in\, H\, \;\text{and}\; \,g_{\,1},\, g_{\,2},\, \,\cdots,\, g_{\,n} \,\in\, K$, where the \,$n$-norms \,$\left\|\,\cdot,\, \cdots,\, \cdot \,\right\|_{1}$\, and \,$\left\|\,\cdot,\, \cdots,\, \cdot \,\right\|_{2}$\, are generated by \,$\left<\,\cdot,\, \cdot \,|\, \cdot,\, \cdots,\, \cdot\,\right>_{1}$\, and \,$\left<\,\cdot,\, \cdot \,|\, \cdot,\, \cdots,\, \cdot\,\right>_{2}$, respectively.\\

The space \,$H_{F} \,\oplus\, K_{G}$\, is the Hilbert space with respect to the inner product:
\[\left<\,p \,\oplus\, q,\, p^{\,\prime} \,\oplus\, q^{\,\prime}\,\right> \,=\, \left<\,p,\, p^{\,\prime}\,\right>_{F} \,+\, \left<\,q,\, q^{\,\prime}\,\right>_{G},\]
for all \,$p,\, p^{\,\prime} \,\in\, H_{F}\; \;\text{and}\; \;q,\, q^{\,\prime} \,\in\, K_{G}$.    

\begin{definition}
Let \,$T \,\in\, \mathcal{B}\left(\,H_{F}\,\right)$\, and \,$U \,\in\, \mathcal{B}\left(\,K_{G}\,\right)$.\,Then the direct sum of the operators \,$T$\, and \,$U$\, is the operator \,$T \,\oplus\, U :\, H_{F} \,\oplus\, K_{G} \,\to\, H_{F} \,\oplus\, K_{G}$\, defined as \,$(\,T \,\oplus\, U\,)\,(\,f \,\oplus\, g\,) \,=\, T\,f \,\oplus\, U\,g$. 
\end{definition}

It is easy to verify that \,$T \,\oplus\, U$\, is a well-defined bounded linear operator whose norm is given by \,$\left\|\,T \,\oplus\, U\,\right\| \,=\, \sup\,\{\,\|\,T\,\|,\, \|\,U\,\|\,\}$.\\

In the following theorem we will show that direct sum of controlled frames is a controlled frame in \,$n$-Hilbert space under some sufficient conditions.

\begin{theorem}
Let \,$\{\,f_{\,i}\,\}_{i \,=\,1}^{\infty}$\, and \,$\{\,g_{\,i}\,\}_{i \,=\,1}^{\infty}$\, be a \,$C_{1}$-controlled and \,$C_{2}$-controlled frames associated to \,$\left(\,a_{\,2},\, \cdots,\, a_{\,n}\,\right)$\, and \,$\left(\,b_{\,2},\, \cdots,\, b_{\,n}\,\right)$\, for \,$H$\, and \,$K$\, with bounds \,$A,\,B$\, and \,$C,\,D$, respectively.\,Then \,$\left\{\,f_{\,i} \,\oplus\, g_{\,i}\,\right\}^{\,\infty}_{i \,=\, 1}$\, is a \,$\left(\,C_{1} \,\oplus\, C_{2}\,\right)$-controlled frame associated to \,$\left(\,a_{\,2} \,\oplus\, b_{\,2},\, \,\cdots,\, a_{\,n} \,\oplus\, b_{\,n}\,\right)$\, for \,$H \,\oplus\, K$, provided for each \,$f \,\in\, H_{F}$\, and \,$g \,\in\, K_{G}$, we have 
\begin{itemize}
\item[$(i)$]$\sum\limits_{i \,=\, 1}^{\,\infty}\,\left<\,f,\, f_{\,i} \,|\, a_{\,2},\, \cdots,\, a_{\,n}\,\right>_{1}\,\left<\,C_{2}\,g_{\,i},\, g \,|\, b_{\,2},\, \cdots,\, b_{\,n}\,\right>_{2} \,=\, 0$, and
\item[$(ii)$]$\sum\limits_{i \,=\, 1}^{\,\infty}\,\left<\,g,\, g_{\,i} \,|\, b_{\,2},\, \cdots,\, b_{\,n}\,\right>_{2}\,\left<\,C_{1}\,f_{\,i},\, f \,|\, a_{\,2},\, \cdots,\, a_{\,n}\,\right>_{1} \,=\, 0$.
\end{itemize}   
\end{theorem}

\begin{proof}
For each \,$f \,\oplus\, g \,\in\, H \,\oplus\, K$, we have
\begin{align*}
&\sum\limits_{i \,=\, 1}^{\,\infty}\,\left<\,f \,\oplus\, g,\, f_{\,i} \,\oplus\, g_{\,i} \,|\, a_{\,2} \,\oplus\, b_{\,2},\, \cdots,\, a_{\,n} \,\oplus\, b_{\,n}\,\right>\,\times\\
&\hspace{1.5cm}\left<\,\left(\,C_{1} \,\oplus\, C_{2}\,\right)\,\left(\,f_{\,i} \,\oplus\, g_{\,i}\,\right),\, f \,\oplus\, g \,|\, a_{\,2} \,\oplus\, b_{\,2},\, \cdots,\, a_{\,n} \,\oplus\, b_{\,n}\,\right>\\
&=\,\sum\limits_{i \,=\, 1}^{\,\infty}\,\left<\,f \,\oplus\, g,\, f_{\,i} \,\oplus\, g_{\,i} \,|\, a_{\,2} \,\oplus\, b_{\,2},\, \cdots,\, a_{\,n} \,\oplus\, b_{\,n}\,\right>\,\times\\
&\hspace{1.5cm}\left<\,C_{1}\,f_{\,i} \,\oplus\, C_{2}\,g_{\,i},\, f \,\oplus\, g \,|\, a_{\,2} \,\oplus\, b_{\,2},\, \cdots,\, a_{\,n} \,\oplus\, b_{\,n}\,\right>\\
&=\,\sum\limits_{i \,=\, 1}^{\,\infty}\,\left\{\,\left<\,f,\, f_{\,i} \,|\, a_{\,2},\, \cdots,\, a_{\,n}\,\right>_{1} \,+\, \left<\,g,\, g_{\,i} \,|\, b_{\,2},\, \cdots,\, b_{\,n}\,\right>_{2}\,\right\}\,\times\\
&\hspace{1.5cm}\left\{\,\left<\,C_{1}\,f_{\,i},\, f \,|\, a_{\,2},\, \cdots,\, a_{\,n}\,\right>_{1} \,+\, \left<\,C_{2}\,g_{\,i},\, g \,|\, b_{\,2},\, \cdots,\, b_{\,n}\,\right>_{2}\,\right\}\\
&=\sum\limits_{i \,=\, 1}^{\,\infty}\,\left<\,f,\, f_{\,i} \,|\, a_{\,2},\, \cdots,\, a_{\,n}\,\right>_{1}\,\left<\,C_{1}\,f_{\,i},\, f \,|\, a_{\,2},\, \cdots,\, a_{\,n}\,\right>_{1} \,+\\
&\hspace{1.5cm}\,\sum\limits_{i \,=\, 1}^{\,\infty}\,\left<\,g,\, g_{\,i} \,|\, b_{\,2},\, \cdots,\, b_{\,n}\,\right>_{2}\,\left<\,C_{2}\,g_{\,i},\, g \,|\, b_{\,2},\, \cdots,\, b_{\,n}\,\right>_{2}\\
&\leq\, B\,\left\|\,f,\, a_{\,2},\, \cdots,\, a_{\,n}\,\right\|^{\,2}_{1} \,+\, D\,\left\|\,g,\, b_{\,2},\, \cdots,\, b_{\,n}\,\right\|^{\,2}_{2}\\
&\leq\, \max\,\{\,B,\,D\,\}\,\left\|\,f \,\oplus\, g,\, a_{\,2} \,\oplus\, b_{\,2},\, \cdots,\, a_{\,n} \,\oplus\, b_{\,n}\,\right\|^{\,2}. 
\end{align*}
On the other hand, for \,$f \,\oplus\, g \,\in\, H \,\oplus\, K$, we have
\begin{align*}
&\sum\limits_{i \,=\, 1}^{\,\infty}\,\left<\,f \,\oplus\, g,\, f_{\,i} \,\oplus\, g_{\,i} \,|\, a_{\,2} \,\oplus\, b_{\,2},\, \cdots,\, a_{\,n} \,\oplus\, b_{\,n}\,\right>\,\times\\
&\hspace{1.5cm}\left<\,\left(\,C_{1} \,\oplus\, C_{2}\,\right)\,\left(\,f_{\,i} \,\oplus\, g_{\,i}\,\right),\, f \,\oplus\, g \,|\, a_{\,2} \,\oplus\, b_{\,2},\, \cdots,\, a_{\,n} \,\oplus\, b_{\,n}\,\right>\\
&\geq\,A\,\left\|\,f,\, a_{\,2},\, \cdots,\, a_{\,n}\,\right\|^{\,2}_{1} \,+\, C\,\left\|\,g,\, b_{\,2},\, \cdots,\, b_{\,n}\,\right\|^{\,2}_{2}\\
&\geq\, \min\,\{\,A,\,C\,\}\,\left\|\,f \,\oplus\, g,\, a_{\,2} \,\oplus\, b_{\,2},\, \cdots,\, a_{\,n} \,\oplus\, b_{\,n}\,\right\|^{\,2}. 
\end{align*}
Thus, the sequence \,$\left\{\,f_{\,i} \,\oplus\, g_{\,i}\,\right\}^{\,\infty}_{i \,=\, 1}$\, is a \,$\left(\,C_{1} \,\oplus\, C_{2}\,\right)$-controlled frame associated to \,$\left(\,a_{\,2} \,\oplus\, b_{\,2},\, \,\cdots,\, a_{\,n} \,\oplus\, b_{\,n}\,\right)$\, for \,$H \,\oplus\, K$\, with bounds \,$\max\,\{\,B,\,D\,\}$\, and \,$\min\,\{\,A,\,C\,\}$. This completes the proof.  
\end{proof}


\begin{thebibliography}{0}
 

\bibitem{B}P. Balazs, J. P. Antonie and A. Grybos, \emph{Weighted and controlled frames: mutual relationship and first numerical properties,} Int. J. Wavelets, Multiresolution Info. Proc., 14 (2010), No. 1, 109-132.

\bibitem{I}I. Bogdanova, P. Vandergheynst, J. P. Antoine, L. Jacques and M. Morrvidone,
\emph{Stereographic wavelet frames on the sphere,} Appl. Comput. Harmon. Anal. 16 (2005), 223-252.

\bibitem{Christensen}O. Christensen,
\emph{An introduction to frames and Riesz bases}, Birkhauser (2008).

\bibitem{Daubechies}
I. Daubechies, A. Grossmann, Y. Mayer,
\emph{Painless nonorthogonal expansions}, Journal of Mathematical Physics 27 (5) (1986) 1271-1283.

\bibitem{Diminnie}
C. Diminnie, S. Gahler, A. White, \emph{2-inner product spaces}, Demonstratio Math. 6 (1973) 525-536.

\bibitem{Duffin}R. J. Duffin, A. C. Schaeffer,
\emph{A class of nonharmonic Fourier series}, Trans. Amer. Math. Soc ., 72, (1952), 341-366.

\bibitem{Folland}G. B. Folland,
\emph{A Course in abstract harmonic analysis}, CRC Press BOCA Raton, Florida.

\bibitem{Gabor}D. Gabor,
\emph{Theory of communications}, J. Inst. Elec. Engrg. 93 (1946), 429-457.

\bibitem{Prasenjit}P. Ghosh and T. K. Samanta, 
\emph{Construction of frame relative to $n$-Hilbert space}, Journal of Linear and Topological Analysis, Vol. 10, No. 02, (2021), 117-130. 

\bibitem{GP}P. Ghosh and T. K. Samanta, 
\emph{Frame in tensor product of \,$n$-Hilbert spaces}, Submitted, arXiv: 2101.01938, Accepted in Sahand Communications in Mathematical Analysis.

\bibitem{PK}P. Ghosh and T. K. Samanta, 
\emph{Atomic systems in $n$-Hilbert spaces and their tensor products}, Submitted, arXiv: 2104.01535.

\bibitem{Gunawan}H. Gunawan,
\emph{On n-inner products, n-norm, and the Cauchy-Schwarz inequality}, Sci. Math. Jpn., 55 (2002), 53-60. 

\bibitem{Mashadi}H. Gunawan, Mashadi,
\emph{On n-normed spaces}, Int. J. Math. Math. Sci., 27 (2001), 631-639. 

\bibitem{Kadison}R. V. Kadison and J. R. Ringrose,
\emph{Fundamentals of the theory of operator algebras}, Vol. I, Academic Press, New York 1983.

\bibitem{Kreyzig}E. Kreyzig,
\emph{Introductory Functional Analysis with Applications.}
Wiley, New York (1989).

\bibitem{Misiak}A. Misiak, \emph{n-inner product spaces}, Math. Nachr., 140(1989), 299-319.

\bibitem{S}S. Rabinson,
\emph{Hilbert space and tensor products}, Lecture notes September 8, 1997. 

\bibitem{Upender}G. Upender Reddy, N. Gopal Reddy and B.\,Krishna Reddy, \emph{Frame operator and Hilbert-Schmidt Operator in Tensor Product of Hilbert Spaces}, Journal of Dynamical Systems and Geometric Theories, 7:1, (2009), 61-70.



\end{thebibliography}
\end{document}